\documentclass[12pt]{article}

\usepackage{amssymb,amsmath,amsthm,amsfonts,latexsym}

\begin{document}

\newtheorem{theorem}{Theorem}[section]
\newtheorem{corollary}[theorem]{Corollary}
\newtheorem{lemma}[theorem]{Lemma}
\newtheorem{proposition}[theorem]{Proposition}
\newtheorem{conjecture}[theorem]{Conjecture}
 \newtheorem{definition}[theorem]{Definition}
 \newtheorem{remark}[theorem]{Remark}

\def\RR{{\mathbb R}}
\def\CC{{\mathbb C}}
\def\ZZ{{\mathbb Z}}
\def\MM{{\mathbb M}}
\def\NN{{\mathbb N}}

\def\A{{\cal A}}
\def\B{{\cal B}}
\def\C{{\cal C}}
\def\F{{\cal F}}
\def\H{{\cal H}}
\def\K{{\cal K}}
\def\N{{\cal N}}
\def\O{{\cal O}}
\def\W{{\cal W}}

\title{{\bf Classification of Subsystems for Graded-Local Nets with
Trivial
Superselection Structure} 
}

\author{
SEBASTIANO CARPI$^*$ \\
Dipartimento di Scienze \\
Universit\`a ``G. d'Annunzio'' di Chieti-Pescara \\
Viale Pindaro 87, 
I-65127 Pescara, Italy \\
{} \\
ROBERTO CONTI\\
Mathematisches Institut\\
Friedrich-Alexander Universit\"at Erlangen-N\"urnberg \\
Bismarckstr. 1 1/2,
D-91054 Erlangen, Germany
}

\date{}
\maketitle
\begin{abstract}
We classify Haag-dual Poincar\'e covariant subsystems
$\B\subset \F$ of a graded-local net $\F$ on 4D Minkowski spacetime
which satisfies standard assumptions and
has trivial superselection structure.
The result applies to the canonical field net $\F_\A$ of a
net $\A$ of local observables satisfying natural assumptions.
As a consequence,
provided that it has no nontrivial internal symmetries,
such an observable net $\A$ is generated by (the abstract
versions of) the local energy-momentum tensor density and the observable
local gauge currents
which appear in the algebraic formulation of the quantum Noether theorem.
Moreover, for a net $\A$ of local observables as
above, we also classify the Poincar\'e covariant local extensions
$\B \supset \A$
which preserve the dynamics.
\end{abstract}

\vfill
\thanks{\noindent $^*$ Partially supported by the Italian MIUR and GNAMPA-INDAM.\\
E-mail: carpi@sci.unich.it, conti@mi.uni-erlangen.de .}

\newpage

\section{Introduction}
\label{introduction}
It is a fundamental insight of the algebraic approach to Quantum Field
Theory  that a proper formulation of  relativistic quantum physics should
be based only on local observable quantities,
see e.g. \cite{Ha}.
The corresponding mathematical structure is a net $\A$ of local
observables, namely an inclusion preserving (isotonous) map which to every
open double cone $\O$ in four
dimensional Minkowski spacetime associates a von Neumann algebra $\A(\O)$
(generated by the observables localized $\O$)  acting on a fixed Hilbert
space $\H_0$ (the vacuum Hilbert space of $\A$) and  satisfying 
mathematically
natural and physically motivated assumptions such as isotony, locality,
Poincar\'e covariance, positivity of the energy and Haag-duality 
(Haag-Kastler axioms).
The charge (superselection) structure of the theory is then encoded in the
representation theory of the quasi-local $C^*$-algebra (still denoted 
$\A$)
which is generated by the local von Neumann algebras $\A(\O)$.
 
The problem was then posed \cite{FOGT,LOPS}, whether it is possible to
reconcile such an approach with the more conventional ones based on
the use of unobservable fields and gauge groups.
 
A major breakthrough was then provided by S.Doplicher and J.E. Roberts in
\cite{DoRo90}.

For any given observable net $\O \mapsto \A(\O)$,
the Doplicher-Roberts reconstruction yields
an associated canonical 
field system with gauge symmetry $(\F,\pi,G)$
describing the superselection structure of the net $\A$ corresponding
to charges localizable in bounded regions (DHR sectors); for details see
\cite{DoRo90} where also the case of topological charges which are
localizable in spacelike cones is considered.
Here $\F=\F_\A$ is the complete normal field net of $\A$,
acting on a larger  Hilbert space $\H \supset \H_0$, the representation
$\pi$ is an embedding of $\A$ into $\F \subset B(\H)$
so that ${\A}={\F}^G$
and the gauge group $G \simeq {\rm Aut}_{\A}(\F)$ is a strongly
compact subgroup of the unitary group $U(\H)$
(to simplify the notation
 we drop the symbol $\pi$ when there is no danger of confusion).
 
Actually any (metrizable) compact group may appear as $G$ \cite{DoPi}.
 
\medskip

If every  DHR sector of $\A$ is Poincar\'e covariant (wich is the
situation considered in this paper) then the net $\F$ is also Poincar\'e
covariant with positive energy. In this case $\A$ is example of covariant
subsystem (or subnet) of $\F$. More generally a covariant subsystem $\B$
of $\F$ is an isotonous map that associate to each double cone $\O$ a
von Neumann subalgebra $\B(\O)$ of $\F(\O)$ which is compatible with
the Poincar\'e symmetry and the grading (giving normal commutation
relations) on $\F$. Besides its natural mathematical interest the study
of covariant subsystems appears also to be useful in the
understanding of the possible
role of local quantum fields of definite physical meaning, such as charge
and energy-momentum densities, in the definition of the net of local
observables (see \cite{sienaCC} and the references therein).
 
In a previous paper \cite{CC} we gave a complete classification of
the (Haag-dual) covariant subsystems of a local field net $\F$ satisfying
some standard additional 
assumptions (like the split property and
the Bisognano-Wichmann property) and having trivial superselection
structure in the sense that every  representation of $\F$ satisfying the
selection criterion of Doplicher, Haag and Roberts \cite{LOPS} (DHR
representation) is unitarily equivalent to a multiple of the vacuum
representation.
The structure that emerged in this analysis is very simple: every 
Haag-dual
covariant subsystem of $\F$ is of the form ${\F_1}^H \otimes 1$ for a
suitable tensor product decomposition $\F=\F_1 \otimes \F_2$ and strongly
compact group $H$ of unbroken internal symmetries of $\F_1$.

Under reasonable assumptions for a net of local
observables $\A$ it was also pointed out in \cite{CC} that the above
result is  sufficient to classify the covariant subsystems of $\F=\F_\A$
(and thus those of $\A$) when the latter net is local
since in that case $\F_\A$ has trivial superselection structure as a
consequence of \cite{CDR}.

The main achievement of this paper is the generalization of the
classification results in \cite{CC}
to the case of graded-local nets i.e. obeying to normal commutation
relations at spacelike distances. Apart from the obvious gain in
mathematical generality our work is intended to remove the
{\it ad hoc} assumption on the locality of $\F_\A$ corresponding to the
absence of DHR sectors of Fermi type, i.e. obeying to (para) Fermi
statistics, for the net of local observables $\A$.

Besides its great theoretical value, the Doplicher-Roberts
(re)construction will provide a major
technical tool for our analysis.
As in \cite{CC} we will repeatedly exploit
the possibility of
comparing such constructions
for different subsystems given by
the
functorial properties of the correspondence $\B \to (\F_\B, G_\B)$
discussed in \cite{CDR}.

\medskip
Compared to \cite{CC},
there are two preliminary problems which have to be settled.

One has to give a meaning to the statement that
``${\cal F}$ has trivial superselection structure''
and this is done by requiring 
that $\F$ has no nontrivial DHR representations,
or, equivalently,
that every DHR representation of the Bose
part $\F^b$ of $\F$ is equivalent to a direct sum of irreducibles 
and that $\F^b$ has only two DHR sectors. 
So in particular
${\cal F} = {\cal F}_{{\cal F}^b}$.
 
One has also to make it clear
in which sense it may hold a tensor product decomposition of ${\cal F}$
(as graded local net) and this is done using the standard notion of
product of Fermionic theories.

Having these differences in mind,
the classification results we obtain (see Theorem \ref{classfull} and
Theorem \ref{classall}) are
nothing but the natural reformulation of the
results in \cite{CC} in the more general context, thus showing
that in the graded local case the structure of
Haag-dual Poincar\'e covariant
subsystems can still be described in terms of internal symmetries 
(cf. \cite{Ar}).
However, though the general strategy is very much the same,
some of the proofs are significantly  different
due to technical complications related to the fact that we did not found
an efficient way to adapt to the graded-local situation
some crucial arguments relying on the work of L. Ge and R.V. Kadison
\cite{GeKa}. These differences are
particularly evident in the proofs of Theorems \ref{classfull} and
\ref{irreducibility} and in fact also provide a partially alternative
argument for the validity of the results in \cite{CC}.
The main new technical ingredients come from the theory of nets of
subfactors \cite{LoRe} and the theory of half-sided modular inclusions of
von Neumann algebras supplemented with some ideas of H.J. Borchers
\cite{Borch99,Wie}.

As a natural application of the classification result we provide a
solution to a problem raised by  S. Doplicher 
(see \cite{Dop})
about the possibility of a
net of local observables to be generated by the corresponding canonical 
local implementations of symmetries 
with the characterization in Theorem
\ref{mettere}.

During our study of subsystems we also realized that some of our methods
can be useful to handle the opposite problem as well, namely to classify
the local extensions of a given observable net $\A$. If a local net
$\B\supset \A$ extends a given local net $\A$ satisfying the same
conditions used in our analysis of subsystems, and if this extension
``preserves the dynamics'', then, modulo isomorphisms, $\B=\F_\A^H$ for a
suitable closed subgroup $H$ of the gauge group of $\A$ (Theorem
\ref{extensiontheorem}).

A crucial assumption on which our results depend and that deserves
some comments is the requirement that the net of local observables
has at most countably many DHR sectors, all with finite statistical
dimension.
 
Although the above properties probably 
are still waiting for
a better understanding
they are strongly supported from the experience: no example of DHR sector
with infinite statistics is known for a theory on a four-dimensional
spacetime and the presence of uncountably many DHR sectors can be ruled out
e.g. by the reasonable requirement that the complete field net fulfills
the split property (the situation is drastically different in the case
of conformal nets on the circle \cite{Car03a,Car03b,Fr,LX}). Thus,
at the present state of knowledge, our results appear to be more than
satisfactory.

We refer the reader to standard textbooks on operator algebras
like \cite{StZs,St,KadRing1,KadRing2}
for all unexplained notions and facts freely used throughout the text.

\section{Preliminaries and assumptions}
\label{preliminaries and assumptions}
We follow closely the discussion in \cite{CC},
pointing out the relevant modifications.

We write $\cal P$ for the component of the identity of the Poincar\'e group
and $\tilde{\cal P}$
for its the universal covering.
Elements of $\tilde{\cal P}$ are denoted by pairs $L=(\Lambda,x)$,
where $L$ is an element of the covering of the connected component of the
Lorentz group and $x$ is a spacetime translation.
$\cal P$ acts in the usual fashion on the
four-dimensional Minkowski spacetime ${\MM}_4$
and the action of $\tilde{\cal P}$ on $\MM_4$ factors through $\cal P$
via the natural covering map $\tilde{\cal P} \to {\cal P}$.

The family of all open double cones and (causal) open wedges in $\MM_4$
will be denoted $\cal K$ and $\cal W$, respectively.
If $S$ is any open region in spacetime,
we denote by $S'$ the interior of the causal complement $S^c$ of $S$.

\medskip
Throughout this paper
we consider
a net $\cal F$ over $\cal K$, namely a correspondence
${\cal O} \mapsto {\cal F}({\cal O})$
between open double cones and von Neumann algebras acting on a fixed
separable
(vacuum) Hilbert space ${\cal H}={\cal H}_{\cal F}$.
The following assumptions
have been widely discussed in the literature
and are by now considered more or less standard:

\begin{itemize}
\item[(i)] {\it Isotony.} If $\O_1 \subset \O_2,\;\O_1,\O_2 \in \K$, then
\begin{equation}
\F(\O_1) \subset \F(\O_2).
\end{equation}

\item[(ii)] {\it Graded locality.}
There exists an involutive unitary operator $\kappa$ on $\cal H$
inducing a net automorphism $\alpha_\kappa$ of $\cal F$,
i.e. $\alpha_\kappa({\cal F}({\cal O})) = {\cal F}({\cal O})$
for each ${\cal O} \in {\cal K}$.
Let
${\cal F}^b(\O) = \{F \in {\cal F}(\O) \ | \ \alpha_\kappa(F)=F\}$
and
${\cal F}^f(\O) = \{F \in {\cal F}(\O) \ | \ \alpha_\kappa(F)=-F\}$
be the even (i.e., Bose)
and the odd (i.e., Fermi) part of $\F(\O)$,
respectively, and let
${\cal F}^t$  be the new (isotonous) net ${\cal F}^b + \kappa {\cal F}^f$
over $\K$.
If $\O_1,\O_2 \in \K$ and $\O_1$ is
spacelike separated from $\O_2$ then
\begin{equation}
\F(\O_1)\subset\F^t(\O_2)^\prime \ .
\end{equation}

\item[(iii)] {\it Covariance.} There is a strongly continuous unitary
representation $U$ of $\tilde{\cal P}$
such that, for every $L\in
\tilde{\cal P}$ and every $\O \in \K$, there holds
\begin{equation}
U(L)\F(\O)U(L)^*=\F(L\O).
\end{equation}
The grading and the spacetime symmetries are compatible, that is
$U(\tilde{\cal P})$ commutes with $\kappa$.

\item[(iv)]{\it Existence and uniqueness of the vacuum.}
There exists a unique (up to a phase) unit vector $\Omega \in \H$
which is invariant
under the restriction of $U$ to the one-parameter subgroup
of spacetime translations.
In addition, one has $\kappa \Omega = \Omega$.

\item[(v)]{\it Positivity of the energy.} The joint spectrum of the
generators of the spacetime translations is contained in the closure
$\overline{V}_+$ of the open forward light cone $V_+$.

\item[(vi)] {\it Reeh-Schlieder property.}
The vacuum vector $\Omega$ is cyclic
and separating
for $\F(\O)$ for every $\O \in \K$.

\item[(vii)] {\it Twisted Haag duality.} For every double cone
$\O \in\K$ there holds
\begin{equation}
\F(\O)^\prime = \F^t(\O^\prime),
\end{equation}
where ,
for every
isotonous net $\F$ and
open set $S \subset \MM_{4}$,
$\F(S)$ denote the von Neumann algebra defined by
\begin{equation}
\F(S)=\bigvee_{\O \subset S} \F(\O).
\end{equation}
Equivalently, one  has
$${\cal F}({\cal O}) =
\bigcap_{{\cal O}_1 \subset {\cal O}'} {\cal F}^t({\cal O}_1)'$$
(in short ${\cal F} = {\cal F}^d$,
where ${\cal F}^d$ is the net defined by the r.h.s.)

\item[(viii)] {\it TCP covariance.} There exists an antiunitary operator
$\Theta$ (the TCP operator) on $\H$ such that:
\begin{align*}
\Theta U(\Lambda,x) \Theta^{-1} & = U(\Lambda,-x)
\quad \forall (\Lambda,x)\in \tilde{\cal P};
\\
\Theta \F(\O) \Theta^{-1} & =\F(-\O)
\quad \forall \O \in \K .
\end{align*}

\item[(ix)] {\it Bisognano-Wichmann property.} Let
\[ W_R=\{x\in \MM_4 : x^1>|x^0| \} \]
be the right wedge and let $\Delta$ and $J$ be the modular operator and
the
modular conjugation of the algebra $\F(W_R)$ with respect to $\Omega$,
respectively. 
Then one has:
\begin{eqnarray}
\Delta^{it}=U(\tilde\Lambda(t),0);
\\
J=Z \Theta U(\tilde{R}_1(\pi),0);
\end{eqnarray}
where $\tilde\Lambda(t)$
and $\tilde{R}_1(\theta)$
denote the lifting in $\tilde{\cal P}$ of the
one-parameter groups 
\begin{equation}
\Lambda(t) 
=\begin{pmatrix}
\cosh 2\pi t & -\sinh 2\pi t & 0 & 0 \\
-\sinh 2\pi t & \cosh 2\pi t & 0 & 0 \\
0 & 0 & 1 & 0 \\
0 & 0 & 0 & 1
\end{pmatrix}
\end{equation} 
of 
Lorentz boosts in the $x^1$-direction
and $R(\theta)$ of spatial rotations around the first axis, respectively,
and $Z=(I+i \kappa)/(1+i)$.

\item[(x)] {\it Split property.} Let $\O_1,\O_2 \in\K$ be open double cones
such that the closure of $\O_1$ is contained in $\O_2$ (as usual we write
$\O_1\subset\subset \O_2$). Then there is a type I factor $\N(\O_1,\O_2)$
such that
\begin{equation}
\F(\O_1)\subset\N(\O_1,\O_2)\subset \F(\O_2).
\end{equation}

\end{itemize}

Assumption (ii) says that $\cal F$ is a graded-local
(or ${\mathbb Z}_2$-graded) net;
for $F \in {\cal F}$ define
$F_+ = (F + \alpha_\kappa(F))/2$ and $F_- = (F - \alpha_\kappa(F))/2$
with $F=F_+ + F_-$,
then given $F_i \in {\cal F}({\cal O}_i)$, $i=1,2$
with ${\cal O}_1$ and ${\cal O}_2$ spacelike separated
the following normal (i.e., Bose-Fermi)
commutation relations hold true:
$$
F_{1+} F_{2+} = F_{2+}F_{1+}, \quad
F_{1+} F_{2-} = F_{2-}F_{1+}, \quad
F_{1-} F_{2-} = - F_{2-}F_{1-} \ .
$$
Note that
${\cal F}^b$,
the net formed by all the elements of ${\cal F}$ which are invariant
under the ${\ZZ}_2$-grading,
is a truly local,
while ${\cal F}^t$ is a graded-local net
(under the same $\kappa$).
Clearly ${\cal F}^{tt}={\cal F}$.

Among the consequences of these axioms one has that
$\F$ acts irreducibly on $\H$, $\F(\MM_4)=B(\H)$.
Moreover, $\Omega$ is $U$-invariant,
and the algebras associated with wedge regions are (type $III_1$) factors.
Strictly speaking, 
one can deduce TCP covariance from the other properties,
see \cite[Theorem 2.10]{GL95}.

By the connection between spin and statistics,
\begin{equation}
\label{spinandstatistics}
\kappa=U(-I,0)
\end{equation}
represents a rotation by angle $2\pi$ about any axis,
see \cite[Theorem 2.11]{GL95}.
(Of course,
in the special case where $\kappa=1$,
we are back in the situation of a local (Bose) net as described in
\cite[Section 2]{CC}.)

From twisted Haag duality it follows that
$\F(\O) = \cap_{\O \subset W}\F(W)$,
thus $\F$ corresponds to an AB-system in the sense of \cite{Wic}.
Note that $Z \F(\O) Z^* = \F^t(\O)$,
for every $\O \in \K$.
Also, the Bisognano-Wichmann property
entails twisted wedge duality,
namely $Z \F(W) Z^* (= \F^t (W)) = \F(W')'$,
see \cite[Prop. 2.5]{GL95}.

Note that $\F^b(S)$ (as defined by additivity)
does not necessarily coincide
with $\F(S)^b := \{F \in \F(S) \ | \ \alpha_\kappa(F)=F\}$
for general open sets $S$.

\medskip
\begin{definition}
A covariant subsystem $\cal B$ of $\cal F$ is an isotonous (nontrivial)
net of von Neumann algebras over $\cal K$, such that
\begin{align*}
{\cal B}({\cal O}) & \subset {\cal F}({\cal O}), \\
U(L){\cal B}({\cal O})U(L)^* & = {\cal B}(L{\cal O})
\end{align*}
for every ${\cal O} \in{\cal K}$ and $L \in \tilde{\cal P}$.
\end{definition}

Then we write ${\cal B} \subset {\cal F}$.
For instance, ${\cal F}^b$ is a covariant subsystem of $\cal F$.
Clearly a covariant subsystem ${\cal B} \subset {\cal F}$ is
local if 
and only if
${\cal B}({\cal O}) \subset {\cal F}^b({\cal O})$ 
for every $\O \in \K$.

For any open set $S \subset {\MM}_4$, we also set
$$\B(S) = 
\bigvee_{{\cal O} \subset S} {\cal B}({\cal O}) \ .$$

If $\B_1$ and $\B_2$ are covariant subsystems of $\F$,
we denote by $\B_1 \vee \B_2$ the covariant subsystem of $\F$ determined
by $(\B_1 \vee \B_2)(\O):=\B_1(\O) \vee \B_2(\O)$, $\O \in \K$.

By the relation (\ref{spinandstatistics})
a covariant subsystem $\B \subset \F$
naturally inherits the grading from $\F$,
namely $\kappa {\cal B}({\cal O}) \kappa^* = {\cal B}({\cal O})$.
Accordingly, $\B^b$ will stand for the local net over $\K$
defined by $\B^b(\O) = \B(\O)^b$.
Also, $\B^t \subset \F^t$ is the (isotonous) 
net $\B^b + \kappa \B^f$.

We denote by $\H_\B := \overline{{\cal B}(\MM_4)\Omega}$
the closed cyclic subspace generated by $\B$ acting on $\Omega$,
and by $E_\B$ the corresponding orthogonal projection of $\H$ onto $\H_\B$.
It follows at once that $E_\B$ commutes with $U$, thus with $\kappa$.

We say that a covariant subsystem ${\cal B} \subset {\cal F}$
is {\it Haag-dual} if
$${\cal B}(\O) = \bigcap_{W \in \W, W \supset \O} {\cal B}(W) \ .$$

As an example, $\F^b$ is a Haag-dual subsystem of $\F$.

Note that a Haag-dual subsystem $\cal B$
does not satisfy twisted Haag duality on $\H_\F$ (unless ${\cal B}=\F$)
however it satisfies (twisted) Haag duality on its own
vacuum Hilbert space $\H_\B$
and the latter property in turn characterizes Haag-dual subsystems.

If $\cal B$ is
Haag-dual, then it satisfies all the properties (i)-(x) 
listed above in restriction to $\H_{\cal B}$
with respect to the restricted representation $\hat{U}$ of $\tilde{\cal P}$
(and grading and TCP operators),
as it can be shown essentially by the same arguments given
in the local case, see \cite[Prop. 2.3]{CC}.
We briefly discuss only the split property.
By repeating the argument in the proof of \cite[Prop. 2.3]{CC},
{\it mutatis mutandis},
it suffices to show that $\H_\B$ is separating for
$\B(\O_1) \vee \B^t(\O'_2) \subset \F(\O_1) \vee \F(\O_2)'$
for every pair of double cones with $\overline{\O}_1 \subset \O_2$.
Pick $\O_0 \subset \O'_1 \cap \O_2$.
Then $\F(\O_0)^b \subset (\B(\O_1) \vee \B^t(\O'_2))'$ 
and $\overline{\F(\O_0)^b \H_\B} \supset \overline{\F(\O_0)^b\Omega} = \H^b 
\ (= \{\xi \in \H \ | \ \kappa \xi = \xi\})$.
Consider $\O \in \K$, $\O \subset \O'_2$,
and pick a Fermi unitary $u$ in $\B(\O)$.
(Such a unitary always exists, 
as it can be seen by applying the polar decomposition 
to any Fermi element in $\B(\O)$ 
and then using Borchers property B for $\B^b$. 
The latter property holds as a consequence of the
split property for $\B^b \subset \F^b$, 
inherited from the split property for $\F$.)
Then $\overline{\F(\O_0)^b u \Omega} = u \H^b = \H^f 
\ (= \{\xi \in \H \ | \ \kappa \xi = -\xi\})$.
Hence $\H_\B$ is cyclic for $(\B(\O) \vee \B^t(\O'_2))'$ and we are done.
\footnote{It is perhaps worth to point out
that the possibly stronger assumption (x$'$) that for every
$\O_1$ and $\O_2$ in $\K$ with $\overline{\O}_1 \subset \O_2$
the triple $(\F(\O_1), \F(\O_2), \Omega)$ is a W$^*$-standard split
inclusion in the sense of \cite{WSSI}
(see Sect. 4) is also inherited by all Haag-dual subsystems.}

If $\cal B$ is not Haag-dual it is always possible to consider
the extension ${\cal B}^d$ defined by
${\cal B}^d(\O) = \cap_{W \in \W, W \supset \O} {\cal B}(W)$,
then ${\cal B}^d$ will be a Haag-dual covariant subsystem of $\F$
with $\B^d(W) = \B(W)$ for every wedge $W \in \W$
and $\H_{{\cal B}^d} = \H_{\cal B}$.
Moreover, in restriction to $\H_{\cal B}$,
$\B^d$ is the (twisted) dual net of $\B$,
namely $\widehat{\B^d}({\cal O})=\hat\B^t({\cal O}')'$ holds on $\H_\B$
where we used $\hat{}$ to denote the restriction of $\B$, resp. $\B^d$
to $\H_\B$.
(However, in order to simplify the notation, 
we shall often write $\B$ in place of $\hat\B$,
specifying if necessary when $\B$ acts on $\H$ or $\H_\B$.)

\medskip
For later use we need to recall the notion of {\it tensor product}
of graded local nets.
Given two graded local nets
${\cal F}_1$ on ${\cal H}_1$ and ${\cal F}_2$ on ${\cal H}_2$
with grading involutive unitaries $\kappa_1$ and $\kappa_2$ respectively,
their tensor product ${\cal F}$
on ${\cal H}_1\otimes{\cal H}_2$
is defined\footnote{The net ${\cal F}_1 \otimes {\cal F}_2$
defined by means of the ordinary tensor product
does not satisfy the normal commutation relations.}
by setting
\begin{equation}
\label{tensorF}
{\cal F}({\cal O})
= {\cal F}_1({\cal O}) \otimes {\cal F}_2({\cal O})^b
+ \kappa_1 {\cal F}_1({\cal O})
\otimes {\cal F}_2({\cal O})^f \ .
\end{equation}
Then $\F$ is still a graded local net
with the diagonal grading $\kappa = \kappa_1 \otimes \kappa_2$,
moreover it satisfies twisted duality
if both the ${\cal F}_i$ do
\cite{Cargese}.
Symbolically we write
${\cal F} = {\cal F}_1 \hat\otimes {\cal F}_2$
to stress that we are dealing with the 
graded tensor product.
\footnote{Other equivalent definitions are possible,
obtained e.g. by exchanging the role of the two components,
or also ${\cal F}^t_1 \hat\otimes {\cal F}_2
= ({\cal F}_1\otimes {\cal F}_2)^b
+ (\kappa_1 \otimes I)({\cal F}_1\otimes {\cal F}_2)^f$.}
One has
${\cal F}^b
= ({\cal F}^t_1 \otimes {\cal F}_2)^b$,
${\cal F}^f
= \kappa_1 \otimes I ({\cal F}^t_1 \otimes {\cal F}_2)^f$.
Note that
$({\cal F}_1 \hat\otimes {\cal F}_2)^{
{\ZZ}_2 \times {\ZZ}_2}
= {\cal F}^b_1 \otimes {\cal F}^b_2$.

With our convention
${\cal F}_1$ sits inside ${\cal F}$ as ${\cal F}_1\otimes I$,
however in general 
$\F_2$ does not share this property.
Nevertheless, $\F_2 \ni F_2 \mapsto 1 \otimes F_{2+} + \kappa_1 \otimes
F_{2-}$ is a 
(normal)
representation of $\F_2$,
unitarily equivalent to $F_2 \mapsto I \otimes F_2$.

For instance,
if the (Fermi) field $\psi_1$ (resp. $\psi_2$)
generates ${\cal F}_1$ (resp. ${\cal F}_2$) then
$\psi_1 \otimes I$ and $\kappa_1 \otimes \psi_2$
will generate ${\cal F}=\F_1 \hat\otimes \F_2$.

The graded tensor product 
$\F = \F_1 \hat\otimes \F_2$ is covariant with respect to the representation $U = U_1 \otimes U_2$
and with vacuum vector $\Omega=\Omega_1 \otimes \Omega_2$,
whenever $\F_i$ is covariant with respect to the representation $U_i$
with vacuum vector $\Omega_i$, $i=1,2$.

\medskip
We say that two graded-local nets
$\F_1$ and $\F_2$ 
as above
are unitarily equivalent (or isomorphic) 
and write 
$\F_1 \simeq \F_2$ if there exists a
unitary operator $W: \H_1 \to \H_2$ such that 
$W \F_1(\O) W^* = \F_2(\O)$ for every $\O \in \K$,
$W \kappa_1 W^* = \kappa_2$,
$W U_1(L) W^* = U_2(L), \ L \in \tilde{\mathcal P}$
and $W\Omega_1 = \Omega_2$.

\medskip
For reader's convenience
we also recall some terminology and few facts that will be used
throughout the paper without any further mention.

A representation $\{\pi,\H_\pi\}$ of the quasi-local 
$C^*$-algebra
associated to a local 
(irreducible, Haag-dual)
net, say $\B$, is said to satisfy the DHR selection criterion,
or simply called a DHR representation,
if for every double cone $\O \in \K$
there exists some unitary $V_\O: \H_\pi \to \H_\B$
such that
\begin{equation}
\label{DHRrep}
V_\O \pi(B) V^*_\O = \pi_0(B), \ B \in \B(\O_1), \ \O_1 \subset \O' \ .
\end{equation}
Here, $\pi_0$ denotes the identical 
(vacuum) representation of $\B$ on $\H_\B$.

Unitary equivalence classes of irreducible DHR representations are called 
{\it DHR superselection sectors} or simply DHR sectors.
The statistics of a DHR sector is described by the statistical dimension,
taking values in ${\mathbb N} \cup \{\infty\}$, and a sign $\pm$
describing the Bose-Fermi alternative.

It is well-known, that a representation $\pi$ satisfies the DHR selection 
criterion if and only if it is unitarily equivalent to some
representation of the form $\pi_0 \circ \rho$, where
$\rho$ is a {\it localized} and {\it transportable} endomorphism
of $\B$.
An account about all this matter can be found e.g. in \cite{Ha},
see also \cite{FOGT,LOPS,Ro}.

In passing, we observe that the definition of DHR representation,
as expressed by equation (\ref{DHRrep}),
carries over to 
graded-local nets, 
however the correspondence with localized and transportable endomorphisms is lost.

\medskip
The results discussed in this paper crucially rely on
the analysis in \cite{CDR}, especially Theorem 4.7 therein.
This provides support for
one further assumption,
which plays an important role in the sequel.

\medskip
(A) \quad
{\it Every representation of the local net $\F^b$
satisfying the DHR selection criterion is a
(possibly infinite)
direct sum of irreducible
representations with finite statistics, moreover
an irreducible DHR representation
that is inequivalent to the vacuum
exists only when $\F^b \subsetneq \F$
and then it is unique
(up to unitary equivalence)}.

\medskip
In particular $\F$ itself is the canonical field net of $\F^b$
in the sense of \cite{DoRo90}.

The following proposition is useful to shed more light on the assumption
(A).
\begin{proposition}
For a field net $\F$ as above, the assumption (A) is satisfied if and only if
every DHR representation of $\F$ is a multiple of the vacuum representation.
\end{proposition}

\begin{proof}
Let $\pi $ be a DHR representation of $\F^b$.
Then $\pi$ is unitarily equivalent to a subrepresentation of the restriction of a DHR representation of $\F$ to
$\F^b$, 
see \cite{CDR}, p.275
(the second paragraph following Proposition 4.3).
If we assume that every DHR representation of $\F$
is a multiple of the vacuum representation,
it follows that $\pi$ is
(equivalent to) a direct sum of irreducible representations of $\F^b$ with finite statistics.
If $\F^b \subsetneq \F$ these are parametrized by $\hat{\mathbb Z}_2 \simeq {\mathbb Z}_2$.

Conversely, 
assume that the condition (A) holds
and let $\tilde\pi$ be a DHR representation of $\F$.
Then, restricting $\tilde\pi$ to $\F^b$, 
it is not difficult to check that 
one gets a DHR representation of $\F^b$.
By assumption, this restriction is thus equivalent to a direct sum of irreducible representations with
finite statistics of $\F^b$.
But then, it follows from \cite[Theorem A.6]{CDR} that
$\tilde\pi$ itself is equivalent to a multiple of the vacuum representation of $\F$.
\end{proof}

Starting from an observable algebra,
the Doplicher-Roberts reconstruction theorem
\cite{DoRo90}
supplies
us with
many
examples of field nets satisfying all the structural assumptions
as above, 
cf. Proposition \ref{directsums}.

\section{Classification of subsystems}
\label{results}
Unless otherwise specified,
throughout this section $\B$ denotes a local Haag-dual covariant
subsystem of ${\cal F}$.

Recall that by \cite[Theorem 3.5]{CDR}
an inclusion of local nets ${\cal A}\subset{\cal B}$
satisfying suitable assumptions
induces an inclusion of the canonical field nets
${\cal F}_{\cal A} \subset {\cal F}_{\cal B}$
compatible with the grading
and thus {\it a fortiori} also
${\cal F}^b_{\cal A} \subset {\cal F}^b_{\cal B}$
and
${\cal F}^f_{\cal A} \subset {\cal F}^f_{\cal B}$.

In particular, in our setting,
from ${\cal B} \subset {\cal F}^b$
we get
${\cal B} \subset
{\cal F}_{\cal B} \subset {\cal F}$
acting on ${\cal H}$
and
${\cal B} \subset {\cal F}^b_{\cal B} \subset
{\cal F}^b$
best considered as acting on ${\cal H}^b$. Moreover these
inclusions are compatible with the action of the Poincar\'e group and
thus,
in particular  ${\cal F}_{\cal B}$ is a covariant subsystem of ${\cal F}$,
cf. \cite{CC} p. 96 and \cite[Theorem 2.11]{sienaCC}.

\medskip
As usual, for a
covariant
subsystem ${\cal B} \subset {\cal F}$
we introduce the {\it coset subsystem} defined by
${\cal B}^c({\cal O})
= {\cal B}(\MM_4)' \cap {\cal F}({\cal O})$, $\O \in \K$.
Then ${\cal B}^c \subset {\cal F}$ is a Haag-dual covariant subsystem.
In principle ${\cal B}^c \subset {\cal F}$
could contain a non-trivial odd part,
in which case it is convenient to consider also
${\cal B}^{cb}
\subset {\cal F}^b$.
If ${\cal B}$ is local then ${\cal B}^c$ could be graded local
but if ${\cal B}$ is
truly
graded local then ${\cal B}^c$ has to be local. As in \cite{CC} we say that
${\cal B}$ is {\it full} in ${\cal F}$ if ${\cal B}^c$ is trivial. Note
however that
in \cite{CDR} the expression ``full'' is used in relation to subsystems
with a different meaning.

\medskip
We take a similar route
as in \cite{CC}.

We denote
$\pi_0$ the vacuum representation of ${\cal B}$
on ${\cal H}_{\cal B}$,
$\pi^0$ that of ${\cal F}$ on ${\cal H}$
and $\pi$ the representation of ${\cal B}$ on ${\cal H}$.
$\pi$ satisfies the DHR selection criterion,
hence
$\pi \simeq \pi_0 \circ \rho$
for some localized and transportable $\rho$.
Note that $\pi_0$ is a subrepresentation of $\pi$,
thus ${\rm id} \prec \rho$. We have the following result, cf.
\cite[Proposition 3.2]{CC}.

\begin{theorem}
\label{basictheorem}
In our situation,
all 
DHR sectors of ${\cal B}$ are
covariant with positive energy
and they have finite statistics.
Furthermore there are at most countably many such DHR sectors
and the actual representation of ${\cal B}$ on ${\cal H}$
is a direct sum of them
in which every DHR sector appears with non zero multiplicity.
\end{theorem}

\begin{proof}
Let $\sigma$ be an irreducible localized transportable endomorphism of
${\cal B}$.

Since ${\cal B}$ and ${\cal F}$ are relatively local
we can extend $\sigma$ to a localized transportable endomorphism
$\hat\sigma$ of ${\cal F}$, 
see \cite[Lemma 2.1]{CDR} (and the paragraph preceding it).

Then, by assumption (A), $\hat\sigma$,
considered as a representation of ${\cal F}$,
is normal on ${\cal F}^b$ and thus is normal on ${\cal F}$ 
by \cite[Theorem A.6]{CDR}.
Since  ${\cal F}$ is irreducible in ${\cal H}$ and
each normal representation of ${\rm B}({\cal H})$ is a multiple of the
identical one we find
$\pi^0 \circ \hat\sigma \simeq \oplus_{i \in I} \pi^0$
for some finite or countable index set $I$.
After restriction of both sides to ${\cal B}$,
$\pi \sigma \simeq \oplus_i \pi$.

From now on the same proof as in \cite[Proposition 3.1]{CC} goes through.
\end{proof}

\begin{proposition}
\label{tensor}
The embedding $\tilde\pi$
of ${\cal F}_{\cal B}$ into ${\cal F}$
satisfies $\tilde\pi \simeq \tilde\pi_0 \otimes I$
where $\tilde\pi_0$ is the vacuum representation of
${\cal F}_{\cal B}$
on ${\cal H}_{{\cal F}_{\cal B}}$.
\end{proposition}

\begin{proof}
By our previous result
the actual representation $\tilde\pi$ of ${\cal F}_{\cal B}$
on ${\cal H}$ is normal with respect to the canonical representation
of ${\cal B}$ on ${\cal H}_{{\cal F}_{\cal B}}$ and thus, by \cite{CDR},
normal
with respect to the actual (irreducible) representation of ${\cal F}_{\cal B}$
on ${\cal H}_{{\cal F}_{\cal B}}$. The conclusion now follows as in the proof
of Theorem \ref{basictheorem} with ${\cal F}_{\cal B}$ instead of ${\cal F}$.
\end{proof}

Besides
the split property for ${\cal F}$ implies that
${\cal F}_{{\cal F}^b_{\cal B}}
= {\cal F}_{\cal B}$,
see \cite{CC}.

\medskip

The following theorem will be of crucial importance. We postpone its
lenghty proof to Appendix \ref{irreducibility of wedges subfactors}.

\begin{theorem}
\label{irreducibility}
If ${\cal F}_{\cal B}$ is full in ${\cal F}$ then
${\cal F}^b_{\cal B}(W) ' \cap {\cal F}(W) = {\CC}I$
for every wedge $W$.
\end{theorem}
We are now ready to state our first classification result
\begin{theorem}
\label{classfull}
Let ${\cal B}$ be a
Haag-dual
local covariant subsystem of ${\cal F}$ and let
${\cal F}_{\cal B}$ be full in ${\cal F}$. Then ${\cal F}_{\cal B}={\cal F}$.
In particular if ${\cal B}$ is full, then there
is a compact group $G$ of unbroken internal symmetries of ${\cal F}$
(with $k \in {\mathcal Z}(G)$, the center of $G$) such that ${\cal B}=
{\cal F}^G$.
\end{theorem}

\begin{proof}
Let us denote ${\cal M}$ the covariant subsystem
${\cal F}^b_{\cal B}$.
By Theorem \ref{irreducibility} we have, for any wedge $W$,
${\cal M}(W) ' \cap {\cal F}(W) = {\CC}I$. Let  $\pi^m_0$ and  $\pi^m$
denote the vacuum representation and the identical representation
on ${\cal H}$ of ${\cal M}$ respectively.
Then, as shown below one can prove that $\pi^m_0$
appears only once in $\pi^m$. Hence
$\tilde\pi \simeq \tilde{\pi}_0$, namely
also the multiplicity of $\tilde{\pi}_0$ in $\tilde\pi$
is one.
Thus ${\cal H}_{{\cal F}_{\cal B}}={\cal H}$ and since
${\cal F}_{\cal B} \subset {\cal F}$ the conclusion follows (e.g. by
twisted Haag duality).
\end{proof}

Fix a wedge $W$
and consider the set
${\cal I}=\{W+a \ | \ a \in {\RR}^4\}$ ordered under inclusion.
Then ${\cal M}(W+a) \subset {\cal F}(W+a)$
defines a directed standard net of subfactors with a standard
conditional expectation as defined in \cite[Section 3]{LoRe}.
Set
$\check{\cal M} = (\cup_a {\cal M}(W+a))^{-\| \cdot \|}$,
$\check{\cal F} = (\cup_a {\cal F}(W+a))^{-\| \cdot \|}$,
then of course
${\cal M}
\subset \check{\cal M}
\subset {\cal M}(\MM_4) = {\cal M}''$
and
${\cal F}
\subset \check{\cal F}
\subset {\cal F}(\MM_4) = {\rm B}({\cal H})$.

Let
$\check{\pi}^0$ denote the representation of $\check{\cal F}$
on ${\cal H}$,
$\check{\pi}^m_0$ that of $\check{\cal M}$ on ${\cal H}_{\cal M}$
and $\check{\pi}^m = \check{\pi}^0|_{\check{\cal M}}$;
clearly
$\pi^0 = \check{\pi}^0 |_{\cal F}$,
$\pi^m_0 = \check{\pi}^m_0 |_{\cal M}$
and $\pi^m = \check{\pi}^m |_{\cal M}$.
By \cite[Corollary 3.3.]{LoRe}
one can construct an endomorphism
$\check\gamma: \check{\cal F} \to \check{\cal M}$
such that $\check{\gamma}|_{{\cal F}(W+a)}$
is Longo's canonical endomorphism of ${\cal F}(W + a)$
into ${\cal M}(W + a)$ whenever $W \subset W+a$,
moreover $\check{\gamma}$ acts trivially on
$\check{\cal M} \cap {\cal F}(W)'$.
It follows from \cite[Proposition 3.4]{LoRe} that
$$\check{\pi}^m \simeq \check{\pi}^m_0 \circ \check{\rho}$$
where $\check\rho = \check\gamma|_{\check{\cal M}}$.
Note that
$\check\rho({\cal M}(W + a)) \subset {\cal M}(W + a)$ if
$W\subset W+a$.

We set $\check\rho_W := \check\rho|_{{\cal M}(W)}$.
It follows from Theorem \ref{irreducibility} and
\cite[Corollary 4.2.]{FiIs} (see also \cite[Theorem 5.2.]{FiIs2}) that
$\iota_{{\cal M}(W)} \prec \check\rho_W$ with multiplicity one.
Let us assume now that $\pi^m_0 \prec \pi^m$ more than once.
Then $\check{\pi}^m_0 \prec \check{\pi}^m$ more than once.
Let $V_1, V_2 \in (\check{\pi}^m_0,\check{\pi}^m_0 \circ \check{\rho})$
be isometries with orthogonal ranges.
Since the net ${\cal M}$ is relatively  local
with respect to ${\cal F}$,  the $C^*$-algebra ${\cal M}_0(W')$ 
generated by all the ${\cal M}(\O)$ with $\O \subset W'$
is contained
$\check{\cal M} \cap {\cal F}(W)'$ and hence the action of $\check{\rho}$
is trivial on it.
For $i = 1,2$ and every $M' \in {\cal M}_0(W')$
we have
$V_i \check{\pi}^m_0 (M')
= \check{\pi}^m_0 \circ \check{\rho}(M')V_i
= \check{\pi}_0 (M')V_i$
i.e.
$V_i \in \check{\pi}^m_0({\cal M}_0(W'))'
= \check{\pi}^m_0 ({\cal M}(W))$
and thus $V_i = \check{\pi}^m_0(W_i)$ with $W_i \in {\cal M}(W)$,
(i=1,2).
But then it follows that $W_1,W_2 \in (\iota_{{\cal M}(W)}, \check\rho_W)$
(notice that $\check{\pi}^m_0$ is faithful)
and they are isometries with orthogonal ranges, which is a contradiction.

 \medskip

We now turn our attention to
to non full subsystems. We begin with the following

\begin{lemma}
\label{bebc}
Let ${\cal B}$ be a 
(not necessarily local)
covariant subsystem of ${\cal F}$. Then
${\cal B} \vee {\cal B}^c$ is full in ${\cal F}$.
\end{lemma}

\begin{proof}
 Since we have
$({\cal B} \vee {\cal B}^c )(\MM_4 )= {\cal B} (\MM_4 ) \vee
{\cal B}^c (\MM_4)$,
then
$({\cal B} \vee {\cal B}^c ) (\MM_4 )' =
{\cal B} (\MM_4 )' \cap {\cal B}^c (\MM_4 )'.$
Hence,  if a double cone ${\cal O}$ is contained in a wedge $W$, then
$$
({\cal B} \vee {\cal B}^c ) (\MM_4)' \cap {\cal F}({\cal O})=
{\cal B}^c ({\cal O}) \cap {\cal B}^c (\MM_4)' \subset
{\cal B}^c (W) \cap {\cal B}^c (W)'.$$
Hence the conclusion follows because ${\cal B}^c (W)$ is a factor.
\end{proof}

\begin{proposition}
\label{fieldcoset}
If ${\cal B}$ is a local covariant subsystem of ${\cal F}$, 
then one has
${\cal F}_{({\cal B}^c)^b} = {\cal B}^c$
on $\H_\F$.
\end{proposition}

\begin{proof}
Let $\tau$ be a transportable endomorphism of $({\cal B}^c)^b$
say localized in the double cone ${\cal O}$
and $\hat{\tau}$ its functorial extension to
${\cal F}^b \supset ({\cal B}^c)^b$.
Since $\hat{\tau}$ is implemented by a 1-cocycle in
$({\cal B}^c)^b \subset {\cal B}'$
we get $\hat\tau(b) = b, \ b \in {\cal B}$,
hence for the corresponding implementing Hilbert space of isometries
in $\F=\F_{\F^b}$ we have
$$H_{\hat\tau} \subset {\cal B}' \cap {\cal F}({\cal O})={\cal B}^c({\cal
O}).$$
Letting $\tau$ to vary in the set $\Delta_{({\cal B}^c)^b}({\cal O})$
of all transportable morphisms of $(\B^c)^b$ localized in $\O$ 
such Hilbert spaces generate ${\cal F}_{({\cal B}^c)^b}({\cal O})$
and then it follows that ${\cal F}_{({\cal B}^c)^b} \subset {\cal B}^c$.
Moreover it is not difficult to see that also the inclusion
${\cal B}^c\subset {\cal F}_{({\cal B}^c)^b} $ holds
and thus the conclusion follows.
\end{proof}

\begin{proposition}
\label{fieldbebc}
The field net
${\cal F}_{{\cal B} \vee ({\cal B}^c)^b}$ acting on its own
vacuum Hilbert space is canonically isomorphic to
$\widehat{{\cal F}_{\cal B}}\hat\otimes \widehat{{\cal B}^c}$
on ${\cal H}_{{\cal F}_{\cal B}}\otimes {\cal H}_{{\cal B}^c}$
as defined in formula (\ref{tensorF})
via the map 
\begin{equation}
\label{canisom}
FB \mapsto 
(\hat{F} \otimes I)(I \otimes \hat{B}_+ + \hat\kappa \otimes \hat{B}_-) \ ,
\end{equation}
where $F \in {\cal F}_{\cal B} 
\subset \F_{\B \vee (\B^c)^b}$ 
and $B = B_+ + B_- \in {\cal B}^c
\subset \F_{\B \vee (\B^c)^b}$.
\end{proposition}

\begin{proof}
Standard arguments relying on the results
in \cite{takesaki} show that the vacuum is a product state for
${\cal B} \vee ({\cal B}^c)^b$ (cf. \cite[Subsection VI.4]{Borch99}).
It follows that the local net ${\cal B} \vee ({\cal B}^c)^b$ acting on
${\cal H}_{{\cal B} \vee ({\cal B}^c)^b}$ is canonically isomorphic to
${\cal B} \otimes ({\cal B}^c)^b$ acting on
${\cal H}_{\cal B}\otimes {\cal H}_{({\cal B}^c)^b}$.
It follows from (the proof of) Theorem \ref{basictheorem} 
(cf. also Proposition \ref{directsums}) that
every factorial DHR representation of $\B$ is a multiple of an irreducible
DHR
representation with finite statistical dimension (in particular it is a
type I representation).
As a consequence every irreducible DHR representation of
${\cal B} \otimes ({\cal B}^c)^b$ is unitarily equivalent to a tensor
product representation and
therefore by Theorem \ref{basictheorem} and Proposition \ref{fieldcoset} 
it is realized up to equivalence 
as a subrepresentation 
on ${\cal H}_{{\cal F}_{\cal B}}\otimes {\cal H}_{{\cal B}^c}$.
Hence,
$\F_{\B \vee (\B^c)^b}$ being generated by the product of Hilbert spaces
of isometries in $\F_\B$, resp. $\B^c$, one has
$\F_{\B \vee (\B^c)^b}=\F_\B \vee \B^c$ on $\H_{\F_{\B \vee (\B^c)^b}}$ and
the conclusion follows from the uniqueness of the canonical field
net \cite[Theorem 3.6]{DoRo90}
along with formula (\ref{tensorF}).

To give more clues on formula (\ref{canisom})
we include a more detailed argument.
It follows from Proposition \ref{tensor}
(with $\F_1=\F_{\B \vee (\B^c)^b}$ in place of $\F$)
the existence of a unitary $W: \H_{\F_1} \to \H_{\F_\B} \otimes \K$,
for some Hilbert space $\K$, such that
$$W F W^* = \hat{F} \otimes I, \quad F \in \F_\B \subset \F_1 \ .$$
Since $\F_\B$ and $(\B^c)^b$ commute
(by an argument similar to the proof of Proposition \ref{fieldcoset}),
it follows that
$$W B W^* =: I \otimes \tau(B), \quad B \in (\B^c)^b \subset \F_1 \ .$$
Moreover, using the easily proven fact that the grading decomposes as a
tensor product, namely
$W \kappa W^* = \hat\kappa \otimes \hat\kappa^c$ 
for some $\hat\kappa^c$,
that $(\B^c)^f$ commutes with $\F^b_\B$ and that $(\B^c)^f (\B^c)^f
\subset (\B^c)^b$ one also deduces that
$$W B W^* =: \hat\kappa \otimes \tau(B), \quad B \in (\B^c)^f \subset \F_1 
\ .$$
Being the representation of $\B^c$ on $\H_{\F_1}$ a multiple of the 
vacuum representation and $\tau$ irreducible, without loss of
generality 
one can identify $\K$ with $\H_{\B^c}$ and $\tau$ with the vacuum
representation of $\B^c$.
But this is exactly formula (\ref{canisom}) and we are done.
\end{proof}

We are now ready to state the complete classification theorem for local
covariant subsystems of ${\cal F}$.

\begin{theorem}
\label{classall}
Let ${\cal B}$ be a local covariant Haag-dual subsystem of ${\cal F}$. Then
$\F = \F_\B \vee \B^c$ on $\H_\F$ and
the net of inclusions $\O \mapsto {\cal B}(\O) \subset {\cal F}(\O)$ 
on ${\cal H}_{\cal F}$ is canonically isomorphic to 
$\O \mapsto \widehat{{\cal F}_{\cal B}}(\O)^H\otimes I 
\subset \widehat{{\cal F}_{\cal B}}(\O)\hat\otimes \widehat{{\cal B}^c}(\O)$ 
on ${\cal H}_{{\cal F}_{\cal B}}\otimes {\cal H}_{{\cal B}^c}$,
where $H$ is the canonical gauge group of $\B$.
\end{theorem}
\begin{proof}
By Proposition \ref{fieldcoset},
we have the following chain of inclusions on $\H_\F$:
$\B \vee \B^c \subset \F_\B \vee \B^c = \F_\B \vee \F_{(\B^c)^b} 
\subset{\cal F}_{{\cal B} \vee ({\cal B}^c)^b} \subset \F$.
Thus,
by Lemma \ref{bebc}
${\cal F}_{{\cal B} \vee ({\cal B}^c)^b}$ is full in $\F$
and by Theorem \ref{classfull} 
we have ${\cal F}= {\cal F}_{{\cal B} \vee ({\cal B}^c)^b}$.
Hence the
conclusion follows from Proposition \ref{fieldbebc}.
\end{proof}

\begin{remark} 
If ${\cal B}$ is an arbitrary covariant
subsystem of ${\cal F}$ satisfying twisted Haag-duality on its 
vacuum Hilbert space ${\cal H}_{\cal B}$ then we can apply
Theorem \ref{classall} to ${\cal B}^b$. Hence the inclusions
${\cal B}^b \subset {\cal B} \subset {\cal F}_{{\cal B}^b}$ allow us
to use Theorem \ref{classall} to classify all (not necessarily local)
covariant subsystems of ${\cal F}$ satisfying  Haag duality or twisted
Haag duality on their own vacuum Hilbert space.
\end{remark}

\section{Nets generated by local generators of symmetries}

In this section we focus our attention on
nets of observables generated by the local generators of symmetries,
i.e. those arising in the framework of the Quantum Noether Theorem \cite{BDL}.
For some background on these nets
we refer the reader to \cite{roberto1,Bob,CC,sienaCC}
(see also \cite{BDLR} for some related issues).
The problem we are interested in is 
to find structural conditions ensuring that 
a given observable net is generated by such 
local generators of symmetries, see \cite{Dop}
(cf. also \cite{LaSc2}). 

\medskip
We consider an observable net $\A$
satisfying the assumption (i)-(x) of Section 2,
with graded locality (resp. twisted Haag duality)
replaced by locality (resp. Haag duality).
Notice that Borchers' property B for $\A$ follows e.g. from the split property.

We further require $\A$ to have countably many DHR 
superselection sectors, all of which have finite statistical dimension.

Let $\F=\F_\A$ and $G = G_\A$
be the canonical field net and the compact gauge group of $\A$.
Arguing as in \cite[Theorem 4.1]{CC},
one can then show that
the assumption (A) for $\F$
introduced in Section 2 is indeed satisfied,
by virtue of \cite[Theorem 4.7]{CDR}.
We record this result here,
as a slight improvement of \cite[Theorem 4.1]{CC}.
\begin{proposition}
\label{directsums}
Let $\A$ be an isotonous net satisfying Haag duality
and the split property 
on its irreducible vacuum representation.
If $\A$ has at most countably many DHR superselection sectors,
all of which have finite statistical dimension,
then any DHR representation of $\A$ is 
unitarily equivalent to 
a (possibly infinite) direct sum of irreducible ones.
Moreover, the canonical field net $\F_\A$ satisfies the assumption (A)
in Sect. 2.
\end{proposition}

All the other properties (i)-(ix)
for $\F$ 
are also satisfied, cf. the discussion in \cite[Subsect. 4.2]{CC}. 
However, in order to apply the analysis in the previous section 
to the present situation
we have 
also to assume property (x) in Sect. 2 to hold for $\F$ since it is still 
unknown whether the split property for $\A$ implies the split property
for $\F_\A$.

Finally in order to construct the local symmetry implementations as in 
\cite{BDL,LASR2} we need to assume that, 
for each pair of double cones $\O_1, \O_2$ with $\overline{\O}_1 \subset 
\O_2$, the vacuum vector $\Omega$ is cyclic for the von Neumann algebra 
$\F(\O_1)'\cap \F(\O_2)$. As a consequence the triple
$(\F(\O_1), \F(\O_2), \Omega)$ is a W$^*$-standard split inclusion
in the sense of \cite{WSSI}.

This last assumption is clearly redundant if $\F$ is local since in 
this case it follows directly from the Reeh-Schlieder property.
Moreover, as shown in the introduction of \cite{LASR2}, in the graded 
local case it would be  a consequence of the split property and the 
Reeh-Schlieder property for $\F$ if $\F^b$ satisfies additivity
and the time slice axiom as in \cite[Sect.1]{LASR1}.

We denote by $\Psi_\Lambda$
the {\it universal localizing map}
associated with the triple 
\begin{equation}
\Lambda = (\F(\O_1), \F(\O_2), \Omega),
\end{equation}
a ${}^*$-isomorphism of $B(\H)$ onto the canonical interpolating factor of
type $I$ between $\F(\O_1)$ and $\F(\O_2)$, see \cite[Sect.3]{BDL}.

\medskip
Let
$$K := G_{\rm max} \equiv
\{k \in U(\H) \ | \ k \F(\O) k^* = \F(\O) \ \forall \O \in \K,
\  k\Omega = \Omega \}
\supset G$$
be the maximal group of
unbroken (unitary) internal symmetries of $\F$,
which in our setting
is automatically strongly compact 
and commutes with Poincar\'e transformations \cite[Theorem 10.4]{WSSI},
and let ${\mathcal C}$ be the net generated by the local version of the
energy-momentum operator \cite{roberto1}, defined by
$$\C(\O) :=
\left( \bigcup_{\O_1, \O_2 : \overline\O_1 \subset \O_2  \subseteq \O}
\Psi_{(\F(\O_1), \F(\O_2), \Omega)} (U(I,\RR^4)) \right) '' \ .$$

The known properties of the universal localizing map 
\cite{BDL,WSSI} imply that $\C$ is a covariant subsystem of $\F$ such 
that
\begin{equation}
\C(\O) \subset \F(\O)^K \subset \F(\O)^G.
\end{equation}
Moreover, by \cite[Corollary 2.5]{DDFL}, we have 
\begin{equation}
U(I,\RR^4) \subset \C(\MM_4) = {\C^d}(\MM_4) \ .
\end{equation}
From
Theorem \ref{classfull}, it readily follows the following result:
\begin{theorem}
The net $\C$ is a full covariant subsystem of $\F$ such that
\begin{equation}
\label{assumption}
{\mathcal C}^d = \F^K \ .
\end{equation}
\end{theorem}
This equality was already proved in \cite{CC}
in the case where $\F$ is Bosonic.
Then one has also $\F=\F_{\C^d}$ and
$K \simeq {\rm Aut}_{\C^d}(\F)$ \cite[Proposition 4.3]{CDR}.

It follows at once that $\A={\mathcal C}^d$ if and only if $G=K$, 
if and only if $\A$ has no proper Haag-dual subsystem full in $\F$,
cf. \cite[Corollary 4.2]{CC}.

\medskip
In the following we somehow exploit the isomorphism between
the lattice structure of 
subgroups of $K$ and that of 
``intermediate'' subsystems of $\F$.

Actually our arguments rest only on the validity of
the equality
(\ref{assumption}).

\medskip
The 
internal symmetries in $K$ leaving $\A = {\F}^G$
globally invariant are exactly those in
${\rm N}_K(G) :=\{k \in K \ | \ kGk^{-1}=G \}$,
the normalizer of $G$ in $K$.
The internal symmetry group of $\A$ can thus be identified with
${\rm N}_K(G) / G$ \cite[Proposition 3.1]{ext}.

We consider the local extension ${\mathcal C}_G$ of $\mathcal C$
in $\F$
by the local currents associated with $G$,
$\C \subset \C_G
\subset {\F}$,
where
$$\C_G(\O) :=
\left( \bigcup_{\O_1, \O_2 : \overline\O_1 \subset \O_2  \subseteq \O}
\Psi_{(\F(\O_1), \F(\O_2), \Omega)} 
(U(I,\RR^4)'' \vee (G'\cap G'')) \right) '' \ .$$
Then ${\mathcal C}^d_G = {\F}^{\tilde{G}}$
for some $\tilde{G} \subset K$
\cite[Theorem 4.1]{CDR};
furthermore,
${\F}={\F}_{\C^d_G}$ and
$\tilde{G} \simeq {\rm Aut}_{\C^d_G}(\F)$,
again by Theorem \ref{classfull}.
Of course, $k \in K$ belongs to $\tilde{G}$
if and only if for each $\Lambda$ one has
\begin{equation}
k \Psi_\Lambda(X) k^{-1}
= \Psi_\Lambda(k X k^{-1})
= \Psi_\Lambda(X)
\ , X \in G' \cap G'' \ , 
\end{equation}
if and only if
\begin{equation}
\label{centerG''}
k X k^{-1} = X , \
X \in G' \cap G''.
\end{equation}
Notice that now $\Lambda$ has disappeared from our condition.
In particular $G \subset \tilde{G}$, so that
${\mathcal C} \subset {\mathcal C}_G
\subset {\F}^G \subset {\F}$.

It is also direct to check that $\tilde{G} \subset {\rm N}_K(G)$,
i.e.
$\tilde{G}$ leaves ${\F}^G$ globally invariant.
To see this, we observe that
the orthogonal projection $E_G$ 
of $\H = \overline{\F \Omega}$
onto
$\H^G = \overline{{\F}^G \Omega}$ lies in
$({\F}^G)' \cap ({\F}^G)'' = G'' \cap G'
\equiv {\mathcal Z}(G'')$, hence $[\tilde{G},E_G]=0$;
then for any $\psi
\in \H^G$ and $k \in \tilde{G}$
one has 
$E_G k \psi = k \psi$, hence $k h k^{-1} \psi = \psi, \ h \in G$
and the conclusion follows since
$$G = \{ k \in K \ | \ k\psi = \psi, \ \psi \in \H^G\} \ .$$

Therefore if $k \in \tilde{G}$ then ${\rm Ad}(k)$ determines
(continuous) automorphisms of both  $G$ and $G''=\A'$.
Let $\alpha_0$ be the natural homomorphism
${\rm N}_K(G) \to {\rm Aut}(G)$,
with kernel ${\rm C}_K(G)$, the centralizer of $G$ in $K$.

By \cite[Lemma 3.13]{DoRo90} all the representations of the compact 
groups above (e.g. $K$) on $\H$ are
quasi-equivalent to the corresponding  left-regular representations.
Since by Eq. (\ref{centerG''}) $k\in K$ belongs to $\tilde{G}$ if and only 
if its adjoint action is trivial on the center of the von Neumann algebra 
$G''$ we can conclude that $k\in \tilde{G}$ if and only if 
$k\in {\rm N}_K(G)$ and $\alpha_0(k) \in {\rm Aut}_{\hat{G}}(G)$, where
${\rm Aut}_{\hat{G}}(G)$ is the group of the automorphisms of $G$ acting 
trivially on the set $\hat{G}$ of equivalence classes of irreducible 
continuous unitary representations of $G$. Hence we have proven the 
following proposition.

\begin{proposition}
We have
$\tilde{G} = \alpha_0^{-1}({\rm Aut}_{\hat{G}}(G)) \ .$
\end{proposition}

One can push this analysis a little bit further.
Dividing
the sequence
$$
\tilde{G} \stackrel i \to {\rm N}_K(G)
\stackrel {\alpha_0} \to {\rm Aut(G)}$$
by $G$
we get another sequence
$$
\tilde{G}/G \stackrel i \to {\rm N}_K(G)/G
\stackrel {\tilde\alpha_0} \to {\rm Aut}(G)/{\rm Inn}(G) =: {\rm Out}(G)$$
so that
\begin{equation}
\label{HtildesuH}
\tilde{G}/G
= {\tilde\alpha}_0^{-1}({\rm Out}_{\hat G}(G))
\end{equation}
where ${\rm Out}_{\hat G}(G) = {\rm Aut}_{\hat G}(G) / {\rm Inn}(G)$.
Notice that
${\rm Ker}(\tilde\alpha_0) = {\rm C}_K(G) / {\mathcal Z}(G)$.
Therefore
$\tilde{G} = G$
if and only if $\tilde{G}/G = \{1\}$,
if and only if
$${\rm C}_K(G) / {\mathcal Z}(G) = \{1\} \quad \mbox{and} \quad
\tilde\alpha_0({\rm N}_K(G)/G) \cap {\rm Out}_{\hat G}(G) = \{1 \} \ .$$

We are now ready
to summarize the above discussion and draw some conclusion. 

It follows from the equation (\ref{HtildesuH}) that
$\tilde{G} / G$
can be identified with the group of
(unbroken)
internal symmetries of the net $\F^G$
that act trivially on the set of its DHR sectors. To see this let 
$\rho$ be an irreducible DHR endomorphism of $\A$ (with finite 
statistical dimension by our previous assumption) localized in a 
double cone $\O \in \K$. Moreover, let 
$H_\rho \subset \F(\O)$ be the corresponding $G-$invariant Hilbert space
of isometries. $H_\rho$ carries an irreducible unitary representation 
$u_\rho$ of $G$ defined by 
$$u_\rho(g) V:= gVg^{-1}, \; V\in H_\rho, g\in G.$$
Now let $k\in {\rm N}_K(G)$ and let $\alpha := {\rm Ad}k|_\A$ be the 
corresponding internal symmetry of $\A$. Then, $\alpha$ acts
on $\rho$ by 
$$\rho \to \rho_\alpha := \alpha \circ \rho \circ \alpha^{-1}.$$
If $V\in H_\rho$ then $kVk^{-1} \in H_{\rho_\alpha}$
(a unitary transformation).
Moreover, 
for every $g\in G$, 
$$u_{\rho_\alpha}(g) kVk^{-1}= gkVk^{-1}g^{-1}=
ku_\rho(\alpha_0(k)(g)) Vk^{-1}.$$ Hence 
$u_{\rho_\alpha}\simeq u_\rho \circ \alpha_0(k)$ and our claim 
follows since if $\rho_1, \rho_2$ are DHR endomorphisms of $\A$ with 
finite statistical dimension then
$\rho_1\simeq \rho_2$ if and only if $u_{\rho_1}\simeq u_{\rho_2}$, 
see \cite{DoRo90}. 

Thanks to the above identification we get immediately the following 
\begin{theorem}
\label{mettere}
Let $\A = \F^G$ be an observable net satisfying all our standing assumptions.
Then one has 
$$\A = \F^G = \C_G^d$$
if and only if $\A$ has no nontrivial (unbroken) internal
symmetries acting identically on the set of its DHR sectors.
In particular, if $\A$ has no nontrivial internal symmetries then 
the above equality of nets holds true.
\end{theorem}
It has been suggested by R. Haag that the existence of nontrivial
internal symmetries for $\A$ is not compatible with the claim that
``the net of observable algebras defines the theory completely without
need of further specifications'' \cite[Sect. IV.1]{Ha}.

\medskip
We briefly list few special cases
for which the computation of $\tilde G$ is 
particularly straightforward:

\begin{itemize}
\item[(i)]
If $G$ is abelian,
then ${\rm Aut}_{\hat G}(G)$ is trivial so that
$\tilde{G}
= \ker(\alpha_0)
= {\rm C}_K(G)$.
For instance, if $K = O(2)$
and $G=SO(2)$ then $\tilde{G}=G$.

\item[(ii)]
If $G$ has no outer automorphisms,
namely ${\rm Aut}(G) = {\rm Inn}(G)$,
\footnote{In the mathematical literature, the groups for which
${\rm Aut}(G) = {\rm Inn}(G)$
and whose center ${\mathcal Z}(G)$ is trivial are called complete. 
}
then $\tilde{G} = {\rm N}_K(G)$.
(The same conclusion holds true if $G$ satisfies the weaker condition that
${\rm Aut}(G) = {\rm Aut}_{\hat{G}}(G)$.)

\item[(iii)] If $G$ is {\it quasi-complete}, meaning that  
$ {\rm Aut}_{\hat{G}}(G) = {\rm Inn}(G)$, then 
$\tilde{G} = G \cdot {\rm C}_K(G)$.
\end{itemize}

In particular 
we have obtained the following result.

\begin{corollary}
If the gauge group $G$ of $\A$ 
has no outer automorphisms
then 
the following conditions are equivalent:

\begin{itemize}
\item[(1)] $G=N_K(G)$, 
namely
$\A=\F^G$ has no nontrivial internal symmetries,

\item[(2)] ${\F}^G = {\mathcal C}^d_G$ \ .
\end{itemize}
\end{corollary}

A special case of this situation can be obtained by taking 
$\F$ to be the net generated by a hermitian scalar free field and
$G = K = {\mathbb Z}_2$,
see \cite[Subsect. 4.1]{CC}, cf. also \cite{LaSc1,LaSc2,Dav}.

\section{Classification of local extensions}
\label{local extensions}
So far we have been dealing only with the
classification problem for subsystems of a given system.
However to some extent
our methods can be useful to handle the extension problem as well.

We assume that a local theory ${\cal A}$ has been given, 
acting on its own vacuum Hilbert space $\H_\A$.

The goal of this section is to setup a framework for classifying all 
the possible
(local) extensions ${\cal B}\supset{\cal A}$
with some additional properties.

The assumptions on $\A$ and $\B$ should allow one to perform the
Doplicher-Roberts reconstruction procedures 
and have a good control on the way these are related.
We show below how this can be achieved
in the case where 
``the energy content of $\B$ is already contained in $\A$''.

In order to be more precise, 
let us assume 
throughout this section 
that $\A$ satisfies 
the same axioms as in the previous section
also including that
it has at most countably many DHR sectors, 
all with finite statistics.
We shall however not need the split property for $\F_\A$.

\begin{definition}
A local extension of the local covariant net $\A$
is a local net $\B$ 
satisfying irreducibility on its separable vacuum Hilbert space $\H_\B$,
Haag duality,
covariance under a representation $V_\B$
fulfilling the spectrum condition and uniqueness of the vacuum,
and containing a covariant subsystem $\A_1$
such that the corresponding net 
$\hat\A_1$ is isomorphic to $\A$.
\end{definition}
In agreement with our notation,
since $\A$ and $\hat\A_1$ are isomorphic,
we shall naturally identify $\A$ and $\A_1$ and write $\A \subset \B$.

In order to state our result, we need two more assumptions.
The first one is of technical nature.
We require the local extension $\B$ as above to satisfy 
the condition of weak additivity, 
namely 
$$\B(\MM_4) = \bigvee_{x \in \RR^4} \B(\O + x)$$
for every $\O \in \K$.
Then the Reeh-Schlieder property and Borchers property B also hold for $\B$.
The second hypothesis, on the energy content, is that 
\begin{equation}
\label{irredplus}
V_\B(I,\RR^4) \subset \A(\MM_4)
\end{equation}
(on $\H_\B$).
As a consequence, $\A$ is full in $\B$.
The meaning of this assumption is to rule out a number of 
unnecessarily ``large'' extensions 
obtained by tensoring the net $\A$ with 
any other arbitrary (local, covariant) net, 
cf. \cite{Car03b}.

\begin{theorem}
\label{extensiontheorem}
Let $\A$ be an observable net satisfying 
all the standing assumptions 
in Section 4.
Let $\B$ be a local extension of $\A$,
satisfying weak additivity and the condition (\ref{irredplus}).

Then $\B$ is 
an intermediate net between $\A$ and 
its canonical field net 
$\F_\A$,
and in fact $\B = \F_\A^H$, the fixed point net of $\F_\A$ 
under the action of some closed subgroup $H$ of the
gauge group of $\A$.
\end{theorem}

\begin{proof}
Since we assume that $\A$ is covariant and the 
spectrum condition holds,
all the DHR sectors of $\A$ are automatically covariant with positive 
energy, see \cite[Theorem 7.1]{GuLo92}.
It follows that $\F_\A = \F_{\A,c}$
where $\F_{\A,c}$ denotes the covariant field net of 
\cite[Section 6]{DoRo90}.

By assumption, 
it is possible to perform the Doplicher-Roberts construction of the
covariant field 
net $\F_{\B,c}$ of $\B$ 
and $\F_\A$ embeds into $\F_{\B,c}$
as a covariant subsystem
\cite[Theorem 2.11]{sienaCC}.

We will make use of Proposition \ref{directsums}.

It is not difficult to see that the natural representation of $\A$ on 
$\H_{\F_{\B,c}}$ is a direct sum of 
DHR representations
(cf. Sect. 3, \cite[Lemma 4.5]{CDR}
and also \cite[Lemma 3.1]{CC})
and thus, by Proposition \ref{directsums}, 
still decomposable into a direct sum of irreducible DHR representations 
with finite statistics.

Arguing similarly as in Sect. 3,
it follows that the representation of $\F_\A={\cal F}_{\A,c}$ 
(thought of as a subsystem of $\F_{\B,c}$)
on the vacuum Hilbert space $\H_{\F_{\B,c}}$ of $\F_{\B,c}$
is a multiple of the vacuum representation of $\F_\A$,
and thus there is a
unitary $W: \H_{\F_{B,c}} \to \H_{\F_\A}\otimes \K$, 
for a suitable Hilbert space $\K$,
such that
$$W F W^* = \hat{F} \otimes I_\K, \ \quad F \in \F_\A \subset \F_{\B,c} $$
(see e.g. the paragraph preceding Theorem 3.1 and Proposition 3.2)
and moreover 
$W U_{\F_{\B,c}} W^* =  U_{\F_{\A}} \otimes \tilde{U}$
for some unitary representation $\tilde{U}$ of $\tilde{\cal P}$ 
on $\K$
with positive energy,
cf. \cite{CC}, p.96. 

Note that, 
by uniqueness of the vacuum, 
if the multiplicity factor $\K$ is nontrivial 
then $\tilde{U}$ has to be nontrivial as well 
and in fact with a unique 
(up to a phase) unit
vector $\Omega_\K$ which is invariant under translations.
To see this, 
without being too much involved in domain problems
which nevertheless can be solved with the help of the spectral theorem,
we consider the situation where $U_2$ and $\tilde{U}$
are unitary representations of $\RR^4$ satisfying the spectrum condition
and $U_1 = U_2 \otimes \tilde{U}$.
Then 
the corresponding generators satisfy the relation 
$$P_1 = P_2 \otimes I + I \otimes \tilde{P} \ .$$
Let us 
assume that $U_i$ has a unique invariant vector $\Omega_i$,
$i=1,2$ (actually $i=1$ would be enough).
Since $P_1 \Omega_1 = 0$, 
the spectrum condition easily implies that 
$(P_2 \otimes I) \Omega_1 = 0$ and $(I \otimes \tilde{P}) \Omega_1 = 0$.
Therefore it follows from the first equation that 
$\Omega_1 = \Omega_2 \otimes \tilde\Omega$ for some vector $\tilde\Omega$
and from the second 
that $\tilde{P}\tilde\Omega = 0$.
Thus $\tilde{U}$ has an invariant vector as well, and this must be unique.
(In alternative, the same result could have been shown with the help of a
Frobenius reciprocity argument.)

Going back to our situation, one can choose $\Omega_\K$ so that
$W \Omega_{\F_{\B,c}} = \Omega_{\F_\A} \otimes \Omega_\K$.

Now,
by the assumption on the energy-momentum operators,
one must have $U_{\F_A} \otimes \tilde{U} = U_{\F_A} \otimes I$ on 
$W\H_\B \subset W\H_{\F_{\B,c}} = \H_{\F_\A}\otimes \K$.
It follows from the last relation that
$W\H_\B \subset \H_{\F_\A} \otimes \Omega_\K 
= \overline{ \hat\F_\A \otimes I (\Omega_{\F_\A} \otimes \Omega_\K) }
$.
Therefore $W B W^* \subset \hat\F_\A \otimes I$, 
i.e. $\B \subset \F_\A$ on $\H_{\F_{\B,c}}$.

Recalling that $U_{\F_{\B,c}} \in \B(\MM_4)$ on $\H_{\F_{\B,c}}$,
one can finally argue that $\tilde{U}$ is trivial on $\K$ and
this immediately yields 
the conclusion that $\K = \CC$ and
\begin{equation}
\label{equfields}
{\cal F}_\A = {\cal F}_{\B,c} \ .
\end{equation}

Thus 
we are back to the situation $\A \subset \B \subset \F_\A$ discussed in 
\cite[Section 4]{CDR} and whence
$$\B=\F^H_\A $$
for some closed subgroup $H$ of $G=G_\A$,
see e.g. \cite[Theorem 4.1]{CDR}.
Thanks to Proposition \ref{directsums}
now
it follows from
\cite[Corollary 4.8]{CDR} that
$$\F_\A = \F_{\B,c} = \F_\B$$
and all the DHR sectors of ${\cal B}$ are covariant
as well.
\end{proof}

\appendix
\section{Proof of Theorem \ref{irreducibility}}
\label{irreducibility of wedges subfactors}
In this appendix we give a proof of Theorem \ref{irreducibility}.
Our proof relies on various ideas from \cite[Sect. V]{Borch99}.
For related problems see \cite[Section 4]{Kos}.

We first recall a convenient notation for the wedges introduced in
\cite{Borch96a} (cf. also \cite{Borch99}). Let $l_1,l_2$ be linearly
independent lightlike vectors in $\overline{V_+}$. Then the region
\begin{equation}
W[l_1,l_2]:= \{\alpha l_1+\beta l_2 + l^\perp :\alpha>0, \beta<0,
l^\perp\cdot l_i =0, i=1,2  \}
\end{equation}
is a wedge and every wedge in $\W$ is of the form $W[l_1,l_2]+a$ for
suitable $l_1, l_2$ and $a\in \MM_4$. Moreover,
$W[l_1,l_2]'=W[l_2,l_1]$. Now let $\F$ be a net satisfying the assumptions
(i)--(x) 
and assumption (A) on triviality of superselection structure 
in Sect. \ref{preliminaries and assumptions}.
Given a wedge $W[l_1,l_2]$ we shall denote by $\Lambda[l_1,l_2](t)$ the
corresponding one-parameter group of Lorentz transformations
(cf. \cite{Borch96a,GL95}).
With this  notation if $\Delta_{[l_1,l_2]}$ is the modular operator for
$(\F(W[l_1,l_2]),\Omega)$ then $\Delta_{[l_1,l_2]}^{it}=
U(\tilde{\Lambda}[l_1,l_2](t),0)$.

In \cite{Borch96a} Borchers considered the intersection of two wedges
$W[l,l_1],W[l,l_2]$ and made the following observations:
\begin{itemize}
\item[a)] $W[l,l_1]\cap W[l,l_2]$ is a nonempty open set;
\item[b)] $\Lambda{[l,l_1]}(t) (W[l,l_1]\cap W[l,l_2]) \subset
W[l,l_1]\cap W[l,l_2]$ for $t\leq 0$.
\end{itemize}
Moreover one finds $\Lambda{[l,l_1]}(t)l=e^{-2\pi t}l$ and
$\Lambda{[l,l_1]}(t)l_1=e^{2\pi t}l_1$.

Now let $S$ be a subset of $\MM_4$ which is contained in some wedge in
$W$. We define (cf. \cite[Sect.III]{TW98})
\begin{equation}
\F^\sharp (S) :=\bigcap_{W\supset S}\F(W).
\end{equation}
It is easy to see that the map $S\mapsto \F^\sharp (S)$ isotonous
and covariant, namely $\F^\sharp (S_1)\subset \F^\sharp (S_2)$ if
$S_1 \subset S_2$ and $U(L)\F^\sharp (S)U(L)^{-1}=\F^\sharp (LS)$ for
every $L\in \tilde{\cal P}$. Clearly $\F^\sharp (W)=\F(W)$ and
$\F^\sharp (\O)=\F(\O)$ for every $W\in \W$ and every $\O \in \K$
but for other regions (even for intersections of family of wedges) the
equality could fail and in general, if $S$ is open and contained in some
wedge, $\F(S)\subset \F^\sharp (S)$.

Similarly, if $\B$ is a Haag-dual covariant subsystem of $\F$ we
define
\begin{equation}
\B^\sharp (S) :=\bigcap_{W\supset S}\B(W).
\end{equation}
Then, the map $S\mapsto \B^\sharp (S)$ is isotonous and covariant and
$\B^\sharp (S)$ coincides with $\B(S)$ if $S\in \W \cup \K$.

Now, given two wedges $W[l,l_1], W[l,l_2]$, the observation of Borchers a)
and b) recalled above and the Bisognano Wichmann property imply
(cf. \cite[Lemma 2.6]{Borch96a}) that the inclusions of von Neumann
algebras
\begin{equation}
\F^\sharp(W[l,l_1]\cap W[l,l_2])\subset \F(W[l,l_i])\; i=1,2,
\end{equation}
are {\it -half-sided modular inclusions} in the sense of
\cite[Definition II.6.1]{Borch99}. Hence, by a theorem of Wiesbrock, Araki
and Zsido (see \cite{Wie} and \cite[Theorem II.6.2]{Borch99}) there are
strongly continuous one-parameter unitary groups $V_i(t), i=1,2$, with
nonnegative generators, leaving the vacuum vector invariant and such that
\begin{equation}
\label{hsmi1}
V_i(t)\F(W[l,l_i])V_i(-t)\subset \F(W[l,l_i]),\; i=1,2,\;t\geq 0,
\end{equation}
\begin{equation}
\label{hsmi2}
\F^\sharp(W[l,l_1]\cap W[l,l_2]) =V_i(1)\F(W[l,l_i])V_i(-1),\;i=1,2.
\end{equation}
Moreover, for $i=1,2$ and $t\in \RR$, $V_i(t)$ is the limit in the strong
operator topology of the sequence
\begin{equation}
\label{trotter}
\left(\Delta_{[l,l_i]}^{-i\frac{t}{2\pi n }}
\Delta_{[l,l_1,l_2]}^{i\frac{t}{2\pi n }} \right)^n,
\end{equation}
where $\Delta_{[l,l_1,l_2]}$ denotes the modular operator of
$\F^\sharp(W[l,l_1]\cap W[l,l_2])$ corresponding to $\Omega$.
As a consequence, $\F^\sharp(W[l,l_1]\cap W[l,l_2])$ is a a factor
of type $III_1$. Note also that by Borchers' theorem
\cite[Theorem II.9]{Borch92} we have, for $i=1,2$ and $t,s\in \RR$,
\begin{equation}
\label{borchtheorem}
\Delta_{[l,l_i]}^{it}V_i(s)\Delta_{[l,l_i]}^{-it}=
V_i(e^{-2\pi t}s)
\end{equation}

The following lemma motivates the introduction of the algebras
$\F^\sharp(S)$.
\begin{lemma}
\label{lemmatwist}
If $l,l_1,l_2$ are lightlike vectors in the closed forward
lightcone such that $l, l_i$ are linearly independent, then the following
holds
\begin{equation}
\F^\sharp(W[l,l_1]\cap W[l,l_2])'=Z\F((W[l,l_1]\cap W[l,l_2])')Z^* 
\end{equation}
\end{lemma}
\begin{proof} We have
\begin{align*}
\F^\sharp(W[l,l_1]\cap W[l,l_2])' & =
\left(\bigcap_{W\supset W[l,l_1]\cap W[l,l_2]}\F(W)\right)' \\
& = \bigvee_{W\supset W[l,l_1]\cap W[l,l_2]}\F(W)' \\
& = \bigvee_{W\supset W[l,l_1]\cap W[l,l_2]}Z\F(W')Z^* \\
& = Z\left(\bigvee_{W\subset ( W[l,l_1]\cap W[l,l_2])'}\F(W)\right)Z^*\\
& = Z\F((W[l,l_1]\cap W[l,l_2])')Z^*,
\end{align*}
where in the last equality we used the convexity of
$W[l,l_1]\cap W[l,l_2]$ and \cite[Proposition 3.1]{TW97}.
\end{proof}

\begin{proposition}
\label{propgenall0}
If $l,l_1,l_2$ are lightlike vectors in the closed forward
lightcone such that $l, l_i$ are linearly independent and $W\in \W$ then
\begin{equation}
\F(\MM_4)= \F(W) \vee \F(W'),
\end{equation}
\begin{equation}
\F(\MM_4)=\F^\sharp(W[l,l_1]\cap W[l,l_2])\vee
\F((W[l,l_1]\cap W[l,l_2])').
\end{equation}
\end{proposition}
\begin{proof} We only prove the second assertion. The proof of the first
is similar but simpler. By Lemma \ref{lemmatwist} we have
\begin{align*}
& \left( \F^\sharp(W[l,l_1]\cap W[l,l_2])\vee  \F((W[l,l_1]\cap
W[l,l_2])') \right)' \\
& = \F^\sharp(W[l,l_1]\cap W[l,l_2])'\cap
\F((W[l,l_1]\cap W[l,l_2])')'\\
& = Z\left(\F^\sharp(W[l,l_1]\cap W[l,l_2])\cap
\F((W[l,l_1]\cap W[l,l_2])') \right)Z^*.
\end{align*}
Now let
$$F\in \F^\sharp(W[l,l_1]\cap W[l,l_2])\cap\F((W[l,l_1]\cap W[l,l_2])')$$
be even with respect to the grading
 (i.e. $\kappa F \kappa$=F). Then, by graded locality,
$$F\in \F^\sharp(W[l,l_1]\cap W[l,l_2])\cap \F^\sharp(W[l,l_1]\cap
W[l,l_2])'$$
and hence $F$ is a multiple of the identity because
 $\F^\sharp(W[l,l_1]\cap W[l,l_2])$ is a
factor. If
$$F \in \F^\sharp(W[l,l_1]\cap W[l,l_2])\cap
\F((W[l,l_1]\cap W[l,l_2])') $$ is odd then $ZFZ^*=i\kappa F$
and hence $i\kappa F$ commutes with $F^*$. It follows that
$$FF^*= -F^*F $$ which implies $F=0$.
\end{proof}

\begin{corollary}
\label{corgenall0} If $\B$ is a local Haag-dual covariant subsystem of
$\F$,  $l,l_1,l_2$ are lightlike vectors in the closed forward
lightcone such that $l, l_i$ are linearly independent and $W\in \W$ then,
on $\H_\F$, we have
\begin{equation}
\F_\B(\MM_4)= \F_\B(W) \vee \F_\B(W'),
\end{equation}
\begin{equation}
\F_\B(\MM_4)=\F_\B^\sharp(W[l,l_1]\cap W[l,l_2])\vee
\F_\B((W[l,l_1]\cap W[l,l_2])').
\end{equation}
\end{corollary}
\begin{proof} By Prop. \ref{propgenall0}
(applied to $\F_\B$ instead of $\F$) the claimed equalities hold on
$\H_{\F_\B}$ and the conclusion follows from Prop. \ref{tensor}.
\end{proof}

\begin{lemma}
\label{connectedlemma}
The set $(W[l,l_1]\cap W[l,l_2])'$ is path connected.
\end{lemma}
\begin{proof}
Recall that, 
given a set $S$, $S' = (S^c)^o$
is defined as the interior of the spacelike complement $S^c$ of $S$.

In order to simplify the notation,
set $W_1 := W[l,l_1]$ and $W_2 := W[l,l_2]$.
Then 
$W_1 \cap W_2 \neq \emptyset$ and
$W'_1 \cap W'_2 \neq \emptyset$. 

One has $(W'_1 \cup W'_2)^c = \overline{W}_1 \cap \overline{W}_2$
\cite[Prop. 2.1, b)]{TW97}
and moreover  
$(W'_1 \cup W'_2)^{cc} 
= ( \overline{W}_1 \cap \overline{W}_2)^c$ is open 
\cite[Prop. 5.6, a)]{TW97}.

We claim that one also has 
$( \overline{W}_1 \cap \overline{W}_2)^c 
= [(W_1 \cap W_2)^c]^o 
\equiv (W_1 \cap W_2)'$.
The inclusion $``\subseteq ''$ is clear.
The opposite inclusion is a consequence of the following three facts:

1) $( \overline{W}_1 \cap \overline{W}_2)^c 
= \cup \{W \ : \ \overline{W}_1 \cap \overline{W}_2 \subset W^c \}$,
as follows by \cite[Theor. 3.2, a)]{TW97}, 
by taking into account the fact that the l.h.s. is actually open;

2) $[(W_1 \cap W_2)^c]^o 
= \cup \{W \ : \ W_1 \cap W_2 \subset W^c\}
= \cup \{W \ : \ \overline{W_1 \cap W_2} \subset W^c\}$,
where 
we use the fact that the spacelike complement of an open set is closed
and the inclusion
$[(W_1 \cap W_2)^c]^o \subset  \cup \{W \ : \ W_1 \cap W_2 
\subset W^c\}$
can be proven with the help of \cite[Prop. 3.1]{TW97};

3) $\overline{W_1 \cap W_2} = (\overline{W}_1 \cap \overline{W}_2)$,
as it can be easily shown recalling that $W_1$ and $W_2$ are open and convex.

Finally, since $W'_1 \cup W'_2$ is connected
it is not difficult to see that
$(W'_1 \cup W'_2)^{cc}=(W_1 \cap W_2)'$
has to be 
connected too.
In fact, if $S$ 
is open and connected and $p \in S^{cc} \setminus S$,
$p$ being spacelike to $S^c$, 
one has that $p$ belongs to the complement of $S^c$.
Since $S$ is open,
this means that the open lightcone
pointed at $p$ intersects $S$ in at least one point, say $q$,
and there is a timelike segment
joining $p$ and $q$ in $S^{cc}$
(cf. the paragraphs preceding \cite[Proposition 2.2]{TW97}).

The proof is complete.
\end{proof}

\begin{proposition}
\label{propexpectation}
 Let $\B$ be a (not necessarily local) Haag-dual covariant subsystem of
$\F$ and let $W[l,l_i], i=1,2$  as above. Then there is a vacuum
preserving normal conditional expectation of
$\F^\sharp(W[l,l_1]\cap W[l,l_2])$ onto
$\B^\sharp(W[l,l_1]\cap W[l,l_2])$.
\end{proposition}
\begin{proof} 
It follows from the Bisognano-Wichmann property and the covariance of $\B$ 
that for every $W\in \W$ the algebra $\B(W)$ is left globally invariant 
by the modular group of $\F(W)$ associated to $\Omega$. Hence, 
by \cite{takesaki}, there is a a vacuum preserving conditional 
expectation $\varepsilon_W$
of $\F(W)$ onto $\B(W)$. Now let $F\in \F^\sharp(W[l,l_1]\cap W[l,l_2])$ 
and $W_a, W_b$ be two wedges containing $W[l,l_1]\cap W[l,l_2]$. Moreover,
let $x_a \in W_a'$ and $x_b \in W_b'$. Since 
$(W[l,l_1]\cap W[l,l_2])'$ is path connected by Lemma 
\ref{connectedlemma} (and open by definition), we can find double cones 
$\O_1,...,\O_n$ all contained in $(W[l,l_1]\cap W[l,l_2])'$ such that 
$x_a\in \O_1$, $x_b \in \O_n$ and $\O_i \cap \O_{i+1} \neq \emptyset$, for 
$i=1,..,n-1$. Then, recalling that $W[l,l_1]\cap W[l,l_2]$ is convex, we 
can use \cite[Prop. 3.1]{TW97} to infer the existence of wedges 
$W_1,...,W_n$ containing $W[l,l_1]\cap W[l,l_2]$ and such that 
$W_1=W_a$, $W_n=W_b$ and $\O_i\subset W_i'$, $i=1,...,n$. Thus, in 
particular, we have $W_i'\cap W_{i+1}' \neq \emptyset$, for 
$i=1,...,n-1$ and hence $\Omega$ is cyclic for 
$Z(\F(W_i')\cap \F(W_{i+1}'))Z^*$ and separating for 
$\F(W_i)\vee \F(W_{i+1})$, $i=1,...,n-1$. Thus, it follows 
from 
$$\varepsilon_{W_i}(F)\Omega =E_\B F\Omega 
=\varepsilon_{W_{i+1}}(F)\Omega,
\; i=1,...,n-1,$$
that $\varepsilon_{W_a}(F)=\varepsilon_{W_b}(F)$. As a consequence the 
restriction of $\varepsilon_W$ to the algebra 
$\F^\sharp(W[l,l_1]\cap W[l,l_2])$
does not depend on $W \supset W[l,l_1]\cap W[l,l_2]$ and gives the 
claimed conditional expectation. 
\end{proof}

As a consequence of the Bisognano-Wichmann property, of Prop.
\ref{propexpectation} and of the results in \cite{takesaki}, for every
Haag-dual covariant subsystem $\B$ of $\F$,
the one-parameter groups $\Delta_{[l,l_i]}^{it}$, $i=1,2$ and
$\Delta_{[l,l_1,l_2]}^{it}$ commute with $E_\B$. Hence, by
(\ref{trotter}) $E_\B$  commutes with $V_i(t)$, $i=1,2$. Moreover,
for $i=1,2$, the following hold
\begin{eqnarray}
\label{EB1}
\B(W[l,l_i])=\F(W[l,l_i])\cap \{E_\B \}', \\
\label{EB2}
\B^\sharp(W[l,l_1]\cap W[l,l_2])=\F^\sharp(W[l,l_1]\cap W[l,l_2])
\cap\{E_\B \}'.
\end{eqnarray}
Hence using Eq. (\ref{hsmi2}) we find
\begin{equation}
\label{hsmiB}
\B^\sharp(W[l,l_1]\cap W[l,l_2]) =V_i(1)\B(W[l,l_i])V_i(-1),\;i=1,2.
\end{equation}

In the following if $M$ is a von Neumann algebra on $\H_\F$ globally
invariant under ${\rm Ad}\kappa(\cdot)$ we shall denote $M^b$ the Bose
part of $M$, i.e. $M^b=M\cap\{ \kappa\}'$. Accordingly if $\B$ is a
covariant subsystem of $\F$ then $\B^b(S)=\B(S)^b$ for every
$S \in \K \cup \W$ and for an arbitrary open set $S$ we have
$\B^b(S)\subset \B(S)^b$.
\begin{proposition}
\label{propisotony}
If $\B$ is a local Haag-dual covariant subsystem of $\F$, $l,l_1,l_2$ are
lightlike vectors in the closed forward lightcone such that $l, l_i$ are
linearly independent and $W\in \W$ then
\begin{equation}
\F_\B(W)'\cap \F(W)^b= \F_\B(\MM_4)'\cap \F(W)^b,
\end{equation}
$$
\F_\B^\sharp(W[l,l_1]\cap W[l,l_2])'\cap
\F^\sharp(W[l,l_1]\cap W[l,l_2])^b$$
\begin{equation}
=\F_\B(\MM_4)'\cap \F^\sharp (W[l,l_1]\cap W[l,l_2])^b.
\end{equation}
\end{proposition}

\begin{proof} We only prove the second
equation. By Corollary \ref{corgenall0} we obtain
$$\F_\B(\MM_4)'\cap \F^\sharp (W[l,l_1]\cap W[l,l_2])^b$$
$$=\F_\B^\sharp(W[l,l_1]\cap W[l,l_2])'\cap
\F^\sharp(W[l,l_1]\cap W[l,l_2])^b \cap
\F_\B((W[l,l_1]\cap W[l,l_2])')'$$
and the conclusion follows since (e.g. by Lemma \ref{lemmatwist})
$$\F^\sharp (W[l,l_1]\cap W[l,l_2])^b
\subset \F_\B((W[l,l_1]\cap W[l,l_2])')'.$$
\end{proof}

\begin{corollary}
\label{corisotony}
 Let $\B$, $l,l_1,l_2$ as in Prop. \ref{propisotony}
and let $W_1, W_2 \in \W$ be such that $W_1\subset W_2$. Then the
following hold
\begin{itemize}
\item[(a)] $\F_\B (W_1)'\cap\F(W_1)^b \subset
\F_\B (W_2)'\cap\F(W_2)^b$,
\item[(b)] $\F_\B^\sharp(W[l,l_1]\cap W[l,l_2])'\cap
\F^\sharp(W[l,l_1]\cap W[l,l_2])^b$ $$ \subset
\F_\B(W[l,l_i])'\cap
\F(W[l,l_i])^b,\;i=1,2.$$
\end{itemize}
\end{corollary}

For a local Haag-dual covariant subsystem $\B$ of $\F$, $l,l_1,l_2$ as in
Prop. \ref{propisotony} and $W\in\W$ we shall now use the following
notations
\begin{equation}
\label{Nnotation1}
N_\B(W):= \F_\B (W)'\cap\F(W)^b,
\end{equation}
\begin{equation}
\label{Nnotation2}
N_\B(W[l,l_1]\cap W[l,l_2]):= \F_\B^\sharp(W[l,l_1]\cap W[l,l_2])'\cap
\F^\sharp(W[l,l_1]\cap W[l,l_2])^b.
\end{equation}

\begin{proposition}
\label{hsmiN}
Let $\B$ be a local Haag-dual covariant subsystem
of $\F$ and let $l,l_1,l_2$ as in Prop. \ref{propisotony}. Then the
following holds

\begin{itemize}
\item[(a)] $V_i(t)N_\B(W[l,l_i])V_i(-t)\subset N_\B(W[l,l_i]),\; i=1,2,
\;t\geq 0,$

\item[(b)] $N_\B(W[l,l_1]\cap W[l,l_2])
=V_i(1)N_\B(W[l,l_i])V_i(-1),\;i=1,2.$
\end{itemize}
\begin{proof} (b) follows easily from Eq. (\ref{hsmi2}) and Eq.
(\ref{hsmiB})
(with $\F_\B$ instead of $\B$). Now, recalling that by (b) in
Corollary \ref{corisotony} we have
$$N_\B(W[l,l_1]\cap W[l,l_2])\subset N_\B(W[l,l_i]),\; i=1,2,$$
it follows from Eq.
(\ref{borchtheorem})
and the covariance of the map $\W \ni W\mapsto N_\B(W)$ that, for every
$s\in\RR$, $i=1,2$,
\begin{align*}
V_i(e^{-2\pi s})N_\B(W[l,l_i])V_i(-e^{-2\pi s}) & =
\Delta_{[l,l_i]}^{is} V_i(1)N_\B(W[l,l_i])V_i(-1)\Delta_{[l,l_i]}^{-is}\\
& = \Delta_{[l,l_i]}^{is}N_\B(W[l,l_1]\cap
W[l,l_2])\Delta_{[l,l_i]}^{-is}\\
& \subset  \Delta_{[l,l_i]}^{is} N_\B(W[l,l_i]) \Delta_{[l,l_i]}^{-is}\\
& = N_\B(W[l,l_i])
\end{align*}
and also (a) is proven.
\end{proof}
\end{proposition}

\begin{lemma}
\label{lemmacoherence1}
Let $\B$ be a local Haag-dual covariant subsystem
of $\F$ and let $l,l_1,l_2$ as in Prop. \ref{propisotony}. Then,
for $i=1,2$, the following holds
\begin{equation}
\overline{N_\B(W[l,l_i])\Omega}=
\overline{N_\B(W[l,l_1]\cap W[l,l_2]) \Omega}.
\end{equation}
In particular
\begin{equation}
\overline{N_\B(W[l,l_1])\Omega}=\overline{N_\B(W[l,l_2])\Omega}.
\end{equation}
\end{lemma}

\begin{proof}
We use a standard Reeh-Schlieder type argument.
By Prop. \ref{hsmiN}, recalling that $V_i(t)\Omega=\Omega$ for
$t\in \RR$, $i=1,2$, we have
\begin{align*}
\overline{N_\B(W[l,l_1]\cap W[l,l_2]) \Omega} & =
V_i(1)\overline{N_\B(W[l,l_i])\Omega} \\
& \subset \overline{N_\B(W[l,l_i])\Omega}.
\end{align*}
Now let $\psi\in (\overline{N_\B(W[l,l_1]\cap W[l,l_2]) \Omega})^\perp$.
Then if $\xi \in \overline{N_\B(W[l,l_i]) \Omega} $
and $t\geq 1$ we have $(\psi,V_i(t)\xi)=0$. Since the self-adjoint
generator of $V_i(t)$ is nonnegative,  the function
$t\mapsto (\psi,V_i(t)\xi)$ is the boundary value of an analytic function
in the upper half-plane. Hence, by the Schwarz reflection principle, it
must vanish for every $t\in \RR$. It follows that $(\psi,\xi)=0$ and hence
that
$ (\overline{N_\B(W[l,l_1]\cap W[l,l_2]) \Omega})^\perp \subset
(\overline{N_\B(W[l,l_i])\Omega})^\perp.$
\end{proof}

\begin{lemma}
\label{lemmacoherence2}
Let $\B$ a local Haag-dual covariant subsystem of $\F$. Then, for every
wedge $W\in \W$, the following holds
\begin{equation}
\overline{N_\B(W)\Omega} =\overline{N_\B(W')\Omega}
\end{equation}
\end{lemma}
\begin{proof}
Let $J_W$ be the modular conjugation of $\F(W)$ with respect to $\Omega$.
Then, by the Bisognano-Wichmann property, $J_W\F(W) J_W =Z\F(W')Z^*$.
Moreover $J_W$ commutes with $\kappa$ and $E_{\F_\B}$. Hence
$$J_W\F_\B(W) J_W =J_W\F(W)\cap\{ E_{\F_\B}\}' J_W = Z\F_\B(W')Z^*$$
and consequently $J_W N_\B(W) J_W =Z N_\B(W')Z^*$.
Accordingly,
$$\overline{N_\B(W)\Omega}  = J_W\overline{N_\B(W)\Omega}
 = Z\overline{N_\B(W')\Omega} = \overline{N_\B(W')\Omega}.$$
\end{proof}

\begin{proposition}
\label{propcoherence}
Let $\B$ a local Haag-dual covariant subsystem of
$\F$. Then, the closed subspace $$\overline{N_\B(W)\Omega}$$ of $\H_\F$
does not depend on the choice of $W\in\W$. Accordingly, the family
$\{N_\B(W): W\in \W \}$ is a coherent family of modular covariant
subalgebras
of $\{\F(W): W\in \W \}$ in the sense of
\cite[Definition VI.3.1]{Borch99}.
\end{proposition}
\begin{proof} Let $l_1,l_2$ and $l'_1, l'_2$ be two pairs of linearly
independent lightlike vectors in the closed forward lightcone. If $l_1$
and $l'_1$ are parallel then $W[l'_1,l'_2] =W[l_1,l'_2]$ and hence
by Lemma \ref{lemmacoherence1}
$$\overline{N_\B(W[l'_1,l'_2] )\Omega}=\overline{N_\B(W[l_1,l_2])\Omega}.$$
On the other hand, if  $l_1$ and $l'_1$ are linearly independent,
using Lemma \ref{lemmacoherence1} and Lemma \ref{lemmacoherence2}, we find
$$\overline{N_\B(W[l_1,l_2] )\Omega}=
\overline{N_\B(W[l_1,l'_1] )\Omega} =
\overline{N_\B(W[l_1,l'_1]' )\Omega}$$
$$=\overline{N_\B(W[l'_1,l_1] )\Omega}=
\overline{N_\B(W[l'_1,l'_2] )\Omega}.$$
To complete the proof it is enough  to show that for any given  wedge
$W\in \W$ the subspace $\overline{N_\B(W+a)\Omega}$ does not
depend on $a\in \MM_4$. To this end we observe that
$\overline{N_\B(W+a)\Omega} = U(I,a) \overline{N_\B(W)\Omega}$
and that, by Corollary \ref{corisotony}, $N_\B(W+a) \subset N_\B(W)$
if $W+a\subset W$. Because of the positivity of the energy, the conclusion
then follows by a Reeh-Schlieder type argument.
\end{proof}
In the following we shall denote the closed subspace $\overline{N_\B(W)\Omega}$
by $\H_{N_\B}$ without any reference to the irrelevant choice of
the wedge $W\in\W$ 
and the corresponding orthogonal projection by $E_{N_\B}$.
We now define for every open double cone $\O\in \K$
\begin{equation}
N_\B(\O) := \bigcap_{W \supset \O} N_\B(W).
\end{equation}
We have the following:
\begin{proposition}
\label{localgeneration}
 For every wedge $W\in \W$ one has
\begin{equation}
N_\B(W) = \bigvee_{\O \subset W} N_\B(\O).
\end{equation}
\end{proposition}
We split the proof of this claim in a series of lemmata.

\begin{lemma}
Let $\O$ be any double cone, then one has
\begin{equation}
N_\B({\O})_{E_{N_\B}}
= \bigcap_{W \supset {\O}} N_\B(W)_{E_{N_\B}} \ .
\end{equation}
\end{lemma}

\begin{proof}
The inclusion ``$\subset$'' is clear.
``$\supset$'': let $X$ denote a generic element in the {\it r.h.s.},
then for every $W \supset {\O}$ there exists $X_W \in N_\B(W)$
such that $X = X_W|_{\H_{N_\B}}$.
We have to show that if $W_a, W_b \supset {\O}$
then $X_{W_a} = X_{W_b}$. Since $\O'$ is connected we can use 
the argument in the proof of Proposition \ref{propexpectation} to find 
wedges $W_1,...,W_n$ containing $\O$, with $W_1=W_a$, 
$W_n=W_b$ and such that $\Omega$ is separating for 
$\F(W_i) \vee \F(W_{i+1})$, $i=1,...,n-1$. The latter property implies 
that $X_{W_i}=X_{W_{i+1}}$, $i=1,...,n-1$, and hence that 
$X_{W_a}=X_{W_b}$.
\end{proof}

\begin{lemma}
For any double cone $\O$
one has $\overline{N_\B(\O)\Omega} = \H_{N_\B}$.
\end{lemma}

\begin{proof}
Consider the family of algebras
$\hat{N}({\cal O}):=N_\B({\cal O})_{E_{N_\B}}$
and
$\hat{N}(W):=N_\B(W)_{E_{N_\B}}$
on $\H_{N_\B}$.
Since
by the previous lemma
$\hat{N}({\O}) = \cap_{W \supset {\O}} \hat{N}(W)$
one deduces that
$$\hat{N}({\O})'
= \bigvee_{W \supset {\O}} \hat{N}(W)'
= \bigvee_{W \subset {\O}'} \hat{N}(W)^t
= (\bigvee_{W \subset {\O}'} \hat{N}(W))^t$$
(in the second equality
we used the fact that $J_W \F(W) J_W = \F^t(W')$).
Now
$\vee_{W \subset \O'} N(W)
\subset \vee_{W \subset \O'} \F(W) = \F(\O')$,
therefore $\Omega$ is separating for $\vee_{W \subset \O'} N(W)$
and henceforth
for $\vee_{W \subset \O'} \hat{N}(W)
= (\vee_{W \subset \O'} N(W))_{E_{N_\B}}$
and also for
$(\vee_{W \subset \O'} \hat{N}(W))^t
= \hat{Z}(\vee_{W \subset \O'}
\hat{N}(W))\hat{Z}^*$.
Thus $\Omega$ is separating for $\hat{N}(\O)'$.
\end{proof}

\begin{lemma}
For every wedge $W$
one has
$\hat{N}(W) = \vee_{\O \subset W} \hat{N}(\O) \ .$
\end{lemma}

\begin{proof}
The inclusion ``$\supset$'' is clear.
Now
let $R$ denote the {\it r.h.s.} and $\sigma$ the modular group of
$(\hat{N}(W),\Omega)$,
then $\sigma_t(R)=R, \ t \in {\mathbb R}$
since $\sigma$ acts like $W$-preserving Lorentz boosts on the double cones
$\O \subset W$.
But $R\Omega$ is dense in $\H_{N_\B}$
and the conclusion follows by Takesaki's theorem on conditional expectations
\cite{takesaki}.
\end{proof}

We are ready to prove the following theorem.
\begin{theorem}
\label{apptheorem}
Let $\B$ be a local Haag-dual covariant subsystem of $\F$
and assume that $\F_\B$ is a full subsystem of $\F$. Then, for every
wedge $W\in \W$, $\F_\B(W)' \cap \F(W)^b =\CC 1$.
\end{theorem}
\begin{proof} First recall that $\F_\B(W)' \cap \F(W)^b =
N_\B(W)$. By Prop. \ref{propisotony} we have
$ N_\B(W) =\F_\B(\MM_4)' \cap \F(W)^b $. Hence, for every $\O\in \K$,
 $$N_\B (\O) \subset \bigcap_{W\supset \O}(\F_\B(\MM_4)' \cap \F(W))
=\F_\B(\MM_4)' \cap \F(\O)=\CC 1$$
and the conclusion follows by Prop. \ref{localgeneration}.
\end{proof}

\begin{corollary}
\label{irreducibility0}
Let $\B$ be a local Haag-dual covariant subsystem of $\F$
and assume that $\F_\B$ is a full subsystem of $\F$. Then, for every
wedge $W\in \W$, $\F_\B(W)' \cap \F(W) =\CC 1$.
\end{corollary}
\begin{proof}
By the Bisognano-Wichmann property the modular group of $\F(W)$ with
respect to $\Omega$ is ergodic (cf. \cite[Lemma 3.2]{longo97}) and leaves
$\F_\B(W)$ globally invariant.
Hence, by \cite{takesaki}, $\F_\B(W)' \cap \F(W)$ has an ergodic modular
group and consequently is either a type III factor or it is trivial but
the first possibility is impossible because of Theorem \ref{apptheorem}.
\end{proof}

To conclude our proof of Theorem \ref{irreducibility} we shall show below
that Corollary \ref{irreducibility0} implies that, if $\B$ be a local
Haag-dual covariant subsystem of $\F$ and  $\F_\B$ is full in $\F$, then
$\F^b_\B(W)' \cap \F(W) =\CC 1$.

 We set $M:=\F(W)$, $N:=\F_\B(W)$ and
denote $\alpha_\kappa$ the automorphism on $M$ induced by the grading operator
$\kappa$. $\alpha_\kappa$ defines an action of $\ZZ_2$ on $M$ and leaves $N$
globally invariant. Let $N_0$ and $M_0$ be the fixpoint algebras for the
action of $\alpha_\kappa$ on $N$ and $M$ respectively. We then have
$\F^b_\B(W)=\F_\B(W)^b=N_0$. Moreover, by Corollary \ref{irreducibility0},
$N\subset M$ is an irreducible inclusion of type III factors and by its
proof $M$ has an ergodic modular group $\sigma^t$ commuting with
$\alpha_\kappa$ and leaving $N$ globally invariant. It follows that also
$N_0 \subset N$ is an irreducible inclusion of type III factors.
By \cite[p. 48]{ILP} $N$ is generated by $N_0$ and a unitary $V$
such that $\alpha_\kappa(V)=-V$. Then $V$ normalizes $N_0$ and $V^2\in N_0$.
Let $\beta := {\rm Ad}V|_{N_0'\cap M}$. Then $\beta$ is an automorphism of
period two. Moreover, the fixpoint algebra $(N_0'\cap M)^\beta$ coincides
with $N'\cap M=\CC1$. Since
$N_0'\cap M$ can be either trivial or a type III factor,
because of the ergodicity of $\sigma^t$,
 we can infer
that $N_0' \cap M= \CC1$, i.e. $\F^b_\B(W)' \cap \F(W) =\CC 1$ and this
concludes the proof of Theorem \ref{irreducibility} .
\bigskip

\medskip
\noindent{\bf Acknowledgements.}
We thank H.-J. Borchers, D. R. Davidson, S. Doplicher, R. Longo,
G. Piacitelli, J. E. Roberts
for several comments and discussions
at different stages of this research.
A part of this work was done while the first named author
(S. C.) was at the Department of Mathematics of the {\it Universit\`a di 
Roma 3}
thanks to a post-doctoral grant of this university.
The final part was carried out while the second named author (R. C.) was
visiting
the Mittag-Leffler Institute in Stockholm during the year devoted to
``Noncommutative Geometry''.
He would like to thank the Organizers for the kind invitation
and the Staff for providing a friendly atmosphere
and perfect working conditions.

\end{document}